\documentclass{article}%
\usepackage[T1]{fontenc}
\usepackage[francais,english]{babel}
\usepackage{amsthm}
\usepackage{amsmath}
\usepackage{amssymb}
\usepackage{mathrsfs}
\usepackage{stmaryrd}
\usepackage{graphicx}
\usepackage[top=3.5cm, bottom=3.5cm, left=3cm, right=3cm, headheight=0.5cm,
footskip=60pt]{geometry}
\usepackage{fancyhdr}
\usepackage{dsfont}
\usepackage{setspace}
\usepackage{listings}
\usepackage{amsfonts}%
\providecommand{\U}[1]{\protect\rule{.1in}{.1in}}
\fancypagestyle{plain}{
	\fancyhead{}
	
}

\newtheorem{theorem}{Theorem}

\newtheorem{definition}{Definition}

\newtheorem{lemma}{Lemma}

\newtheorem{proposition}{Proposition}
\newtheorem{remark}{Remark}

\numberwithin{lemma}{section}
\numberwithin{proposition}{section}
\numberwithin{equation}{section}
\numberwithin{remark}{section}
\numberwithin{definition}{section}
\allowdisplaybreaks
\begin{document}

\title{New long time existence results for a class of Boussinesq-type systems}
\author{Cosmin Burtea\footnote{ Corresponding author. Email address:
cosmin.burtea@univ-paris-est.fr} \footnote{This work was supported by a grant of the Romanian National Authority for
Scientific Research and Innovation, CNCS - UEFISCDI, project number
PN-II-RU-TE-2014-4-0320}\\Universit\'e Paris-Est Cr\'eteil, LAMA - CNRS UMR 8050,\\61 Avenue du G\'{e}n\'{e}ral de Gaulle, 94010 Cr\'{e}teil, France}
\maketitle

\begin{abstract}
In this paper we deal with the long time existence for the Cauchy problem
associated to some asymptotic models for long wave, small amplitude gravity
surface water waves. We generalize some of the results that can be found in
the literature devoted to the study of Boussinesq systems by implementing an
energy method on spectrally localized equations. In particular, we obtain
better results in terms of the regularity level required to solve the initial
value problem on large time scales and we do not make use of the positive
depth assumption.

\begin{description}
\item[Keywords] Boussinesq systems; Long time existence;

\end{description}

\section{Introduction}

\end{abstract}

\subsection{The abcd systems}

The following $abcd$ Boussinesq systems were introduced in \cite{Bona1} as
asymptotic models for studying long wave, small amplitude gravity surface
water waves:%

\begin{equation}
\left\{
\begin{array}
[c]{l}%
\left(  I-\varepsilon b\Delta\right)  \partial_{t}\eta+\operatorname{div}%
V+a\varepsilon\operatorname{div}\Delta V+\varepsilon\operatorname{div}\left(
\eta V\right)  =0,\\
\left(  I-\varepsilon d\Delta\right)  \partial_{t}V+\nabla\eta+c\varepsilon
\nabla\Delta\eta+\varepsilon\frac{1}{2}\nabla\left\vert V\right\vert ^{2}=0.
\end{array}
\right.  \label{eq1}%
\end{equation}
In system \eqref{eq1} $\varepsilon$ is a small parameter while
\[
\left\{
\begin{array}
[c]{c}%
\eta=\eta\left(  t,x\right)  \in\mathbb{R},\\
V=V\left(  t,x\right)  \in\mathbb{R}^{n},
\end{array}
\right.
\]
with $\left(  t,x\right)  \in\left[  0,\infty\right)  \times\mathbb{R}^{n}$
are approximations of the free surface of the water and of the fluid's
velocity respectively. As it will soon be clearer, we mention that the only
values of $n$ for which \eqref{eq1} is physically relevant are $n=1,2$. The
above family of systems is derived from the classical mathematical formulation
of the water waves problem: considering a layer of incompressible,
irrotational, perfect fluid flowing through a canal with flat bottom
represented by the plane:%
\[
\{\left(  x,y,z\right)  :z=-h\},
\]
where $h>0$ and assuming that the free surface resulting from an initial
perturbation of the steady state can be described as being the graph of a
function $\eta$ over the flat bottom, the water waves problem is governed by
the following system of equations:%
\begin{equation}
\left\{
\begin{array}
[c]{rrr}%
\Delta\phi+\partial_{z}^{2}\phi=0 & \text{in} & -h\leq z\leq\eta\left(
x,y,t\right)  ,\\
\partial_{z}\phi=0 & \text{at} & z=-h,\\
\partial_{t}\eta+\nabla\phi\nabla\eta-\partial_{z}\phi=0 & \text{at} &
z=\eta\left(  x,y,t\right)  ,\\
\partial_{t}\phi+\frac{1}{2}\left(  \left\vert \nabla\phi\right\vert
^{2}+\left\vert \partial_{z}\phi\right\vert ^{2}\right)  +gz=0 & \text{at} &
z=\eta\left(  x,y,t\right)
\end{array}
\right.  \tag{$\mathcal{WW}$}\label{WW}%
\end{equation}
where $\phi$ stands for the fluid's velocity potential and $g$ is the
acceleration of gravity. The operators $\Delta$ and $\nabla$ are taken with
respect to $(x,y)$. Of course, in many applications the above system of
equations raises a significant number of problems both theoretically and
numerically. This is the reason why an important number of approximate models
have been established, each of them dealing with some particular physical
regimes. The $abcd$ systems of equations deals with the so called Boussinesq
regime which we will explain now. We consider the following quantities:
$A=\max_{x,y,t}\left\vert \eta\right\vert $ the maximum amplitude encountered
in the wave motion, $l$ the smallest wavelength for which the flow has
significant energy and $c_{0}=\sqrt{gh}$ the kinematic wave velocity. The
Boussinesq regime is characterized by the parameters:%
\begin{equation}
\alpha=\frac{A}{h},\text{ }\beta=\left(  \frac{h}{l}\right)  ^{2},\text{
}S=\frac{\alpha}{\beta}, \label{Bous}%
\end{equation}
which are supposed to obey the following relations:%
\[
\alpha\ll1,\text{ }\beta\ll1\text{ and }S\approx1.
\]
Supposing actually that $S=1$ and choosing $\varepsilon=\alpha=\beta$, the
systems \eqref{eq1} are derived back in \cite{Bona1}. The unknown functions
$\left(  \eta,V\right)  $ in \eqref{eq1} represent the deviation of the free
surface from the rest state while $V$ is an $O\left(  \varepsilon^{2}\right)
$ approximation of the velocity $\nabla\phi$ taken at a certain depth.
Actually, the zeros on the right hand side of \eqref{eq1} represent the
$O\left(  \varepsilon^{2}\right)  $ terms neglected in establishing
\eqref{eq1}. The parameters $a,b,c,d$ are also restricted by:%
\begin{equation}
a+b+c+d=\frac{1}{3}. \label{suma}%
\end{equation}
Asymptotic models taking into account different topographies of the bottom
were also derived, see for instance \cite{Chaz}, Section $2$, for bottoms
given by the surface:%
\[
\{\left(  x,y,z\right)  :z=-h+\varepsilon S(x,y)\},
\]
where $S$ is smooth enough or \cite{Chen} for slowly variating bottoms i.e.
the function $S=S\left(  t,x,y\right)  $ depends also on time. A systematic
study of asymptotic models for the water waves problem along with their
rigorous justification can be found for instance in \cite{Lannes}.

The study of the local well-posedness of the $abcd$ systems is the subject
under investigation in several papers see for instance \cite{Bona1},
\cite{Bona2}, for space dimension $n=1$ or \cite{Anh}, \cite{Bona3} (the
BBM-BBM case $b=d>0$, $a=c=0$), \cite{Saut2} (the KdV-KdV case $b=d=0$,
$a=c>0$) for dimension $n=2$. In \cite{Bona1} it is shown that the linearized
equation near the null solution of \eqref{eq1} is well posed in two generic
cases, namely:%
\begin{align}
a  &  \leq0,\text{ }c\leq0,\text{ }b\geq0,\text{ }d\geq0\label{1111}\\
\text{or }a  &  =c\geq0\text{ and }b\geq0,\text{ }d\geq0. \label{2222}%
\end{align}
It was proved in \cite{Bona4} (see also \cite{Lannes}) that the error estimate
between the solution of \eqref{eq1} and the water wave system is cumulating in
time like $O\left(  \varepsilon^{2}t\right)  $. Thus, solutions of \eqref{eq1}
that exist on time intervals of order $O\left(  \frac{1}{\varepsilon}\right)
$ are good approximation for \eqref{WW} as the error remains of small order
i.e. $O\left(  \varepsilon\right)  $. Actually, due to the previous mentioned
error estimate, on time scales larger then $O\left(  \frac{1}{\varepsilon^{2}%
}\right)  $ the solution $\left(  \eta,V\right)  $ stops being relevant as an
approximation for the original system. The question of long $O\left(  \frac
{1}{\varepsilon}\right)  $-time existence of solutions of \eqref{eq1} did not
receive a satisfactory answer until in \cite{Saut1} where, the case
\eqref{1111} was treated and long time existence for the Cauchy problem was
systematically proved, provided that the initial data lies in some Sobolev
spaces. The difficulty comes from the lack of symmetry (the $\varepsilon
\eta\operatorname{div}V$ term from the first equation of the $abcd$-system) of
\eqref{eq1}. Because of the dispersive operators $-\varepsilon b\Delta
\partial_{t}+a\operatorname{div}\Delta$, $-\varepsilon d\Delta\partial
_{t}+c\nabla\Delta$, classical symmetrizing techniques for hyperbolic systems
of PDE's fail to be successful.

Global existence is known to hold true for \eqref{eq1} in dimension $n=1$ for
some particular cases: the case
\[
a=b=c=0,\text{ }d>0,
\]
studied by Amick \cite{Amick} and Schonbek \cite{Schon} and in the case%
\[
b=d>0,\text{ }a\leq0,~c<0,
\]
assuming some smallness condition on the initial data, see \cite{Bona2}. In
both cases, it is assumed that:%
\[
\inf_{x\in\mathbb{R}}\left\{  1+\varepsilon\eta_{0}\left(  x\right)  \right\}
>0
\]
a condition that makes perfect sense from a physical point of view as
$1+\varepsilon\eta\left(  t,x\right)  $ represents the total height of the
water above the bottom in $x$ at time $t$ . The latter case uses the
Hamiltonian structure of the system, namely, when $b=d$ we have:%
\begin{equation}
\frac{d}{dt}\left(  \int\eta^{2}+\left(  1+\varepsilon\eta\right)
V^{2}-\varepsilon c\left(  \partial_{x}\eta\right)  ^{2}-\varepsilon a\left(
\partial_{x}V\right)  ^{2}\right)  dx=0. \label{Hamil}%
\end{equation}
\ \

\subsection{The main result}

The aim of this paper is to generalize most of the results presented in
\cite{Saut1}, namely, we address the long time existence issue for the general
$abcd$ systems. More precisely we wish to construct solutions for \eqref{eq1}
for which the time of existence is bounded from below by a $O\left(  \frac
{1}{\varepsilon}\right)  $-order quantity. The key ingredient is that we
establish energy-type estimations for spectrally localized equations and by
doing so we require a lower regularity index than\ the one used in
\cite{Saut1}. Also, we avoid using the non-cavitation condition%
\begin{equation}
1+\varepsilon\eta_{0}\left(  x\right)  \geq\alpha>0, \label{ip}%
\end{equation}
imposed on the initial datum $\eta_{0}$ or the curl free condition on the
initial data $V_{0}$ used in \cite{Bona3}. As opposed to the method used in
\cite{Saut1}, ours permits us to treat in a unified manner most of the cases
corresponding to the values of the $a,b,c,d$ parameters. In order to carry on
our approach we need some restrictions on the $a,b,c,d$ parameters. More
precisely, we will consider the case \eqref{1111} such that $b+d>0$ excluding
the following two cases:%
\begin{equation}
\left\{
\begin{array}
[c]{l}%
a=d=0,\text{ }c<0,b>0.\\
a=b=0,\text{ }c<0,d>0.
\end{array}
\right.  \label{abcd}%
\end{equation}
In Section \ref{sub1} we put forward the basic ingredients in order to obtain
long time existence for the former two cases. In view of \eqref{suma}, $b+d>0$
is not restrictive as far as $a,c$ are less or equal to $0$.

First, we will slightly change the form of \eqref{eq1} noticing that the
divergence free part of $V$ remains constant during time. Indeed, formally, if
$\left(  \eta,\bar{V}\right)  $ is a solution of \eqref{eq1} with initial data%
\[
\eta_{|t=0}=\eta_{0}\text{ , }\bar{V}_{|t=0}=\bar{V}_{0},
\]
then%
\begin{equation}
\partial_{t}\bar{V}=-\nabla\left(  I-\varepsilon d\Delta\right)  ^{-1}\left(
\eta+c\varepsilon\Delta\eta+\varepsilon\frac{1}{2}\left\vert \bar
{V}\right\vert ^{2}\right)  , \label{rot}%
\end{equation}
and consequently we have that:%
\[
\partial_{t}\mathbb{P}\bar{V}=0,
\]
where%
\[
\mathbb{P}\bar{V}=\mathcal{F}^{-1}\left(  \left(  I-\sum_{i=1,n}\frac{\xi
_{i}\xi}{\left\vert \xi\right\vert ^{2}}\right)  \mathcal{F}\left(  \bar
{V}\right)  \right)
\]
is the Leray projector over divergence free vector fields. Thus, we get that:%
\[
\mathbb{P}\bar{V}=\mathbb{P}\bar{V}_{0}\overset{not}{=}W.
\]
Of course, we have that%
\[
\operatorname{div}W=0.
\]
By setting%
\[
\left\{
\begin{array}
[c]{l}%
\bar{V}_{0}=V_{0}+W\text{ ,}\\
\bar{V}=V+W
\end{array}
\right.
\]
we infer that the system verified by the couple $(\eta,V)$ is the following:%
\begin{equation}
\left\{
\begin{array}
[c]{l}%
\left(  I-\varepsilon b\Delta\right)  \partial_{t}\eta+\operatorname{div}%
V+a\varepsilon\operatorname{div}\Delta V+\varepsilon W\nabla\eta
+\varepsilon\operatorname{div}\left(  \eta V\right)  =0,\\
\left(  I-\varepsilon d\Delta\right)  \partial_{t}V+\nabla\eta+c\varepsilon
\nabla\Delta\eta+\varepsilon\frac{1}{2}\nabla\left\vert W\right\vert
^{2}+\varepsilon\nabla WV+\varepsilon\nabla VW+\varepsilon\frac{1}{2}%
\nabla\left\vert V\right\vert ^{2}=0,\\
\eta_{|t=0}=\eta_{0}\text{ , }V_{|t=0}=V_{0}.
\end{array}
\right.  \label{eq2}%
\end{equation}
Also because of \eqref{rot} we get that $\operatorname{curl}V=0$ at any time
meaning that for all $l,k\in\overline{1,n}$ we have:%
\begin{equation}
\partial_{l}V^{k}=\partial_{k}V^{l}. \label{curl}%
\end{equation}
The advantage of working with system \eqref{eq2} instead of \eqref{eq1} is
two-folded. On one side certain commutators involving $V$ behave better if its
curl is $0$ and on the other side, for some values of the $a,b,c,d$ parameters
e.g. $a=c=d=0$, $b>0$ we need less regularity on the initial data $V_{0}$.

Let us establish some notations. The spaces $L^{p}\left(  \mathbb{R}%
^{n}\right)  $ with $p\in\lbrack1,\infty]$ will denote the classical Lebesgue
spaces. Given $s\in\mathbb{R}$ we will consider the following set of indices:%

\begin{equation}
\left\{
\begin{array}
[c]{l}%
s_{1}=s+sgn\left(  b\right)  -sgn\left(  c\right)  ,\\
s_{2}=s+sgn\left(  d\right)  -sgn\left(  a\right)  ,\\
s_{3}=s+1-sgn\left(  a\right)  ,
\end{array}
\right.  \label{relatie}%
\end{equation}
where the sign function $sgn$ is given by:%
\[
sgn\left(  x\right)  =\left\{
\begin{array}
[c]{rrr}%
1 & \text{if} & x>0,\\
0 & \text{if} & x=0,\\
-1 & \text{if} & x<0,
\end{array}
\right.
\]
and $a,b,c,d$ are chosen as in \eqref{1111}. We will denote by $H^{s}\left(
\mathbb{R}^{n}\right)  $ the classical Sobolev space of regularity index $s$
with the norm%
\begin{equation}
\left\Vert \eta\right\Vert _{H^{s}}^{2}=%
{\displaystyle\int\limits_{\mathbb{R}^{n}}}
\left(  1+\left\vert \xi\right\vert ^{2}\right)  ^{s}\left\vert \hat{\eta
}\left(  \xi\right)  \right\vert ^{2}d\xi. \label{hs}%
\end{equation}
For any vector, matrix or $3$-tensor field with components in $H^{s}\left(
\mathbb{R}^{n}\right)  $ we denote its Sobolev norm by the square root of the
sum of the squares of the Sobolev norms of its components. For any pair
$\left(  \eta,V\right)  \in H^{s_{1}}\left(  \mathbb{R}^{n}\right)
\times\left(  H^{s_{2}}\left(  \mathbb{R}^{n}\right)  \right)  ^{n}$ we will
use the notation%
\begin{align*}
\left\Vert \left(  \eta,V\right)  \right\Vert _{s}^{2}  &  =\left\Vert
\eta\right\Vert _{H^{s}}^{2}+\varepsilon(b-c)\left\Vert \nabla\eta\right\Vert
_{H^{s}}^{2}+\varepsilon^{2}(-c)b\left\Vert \nabla^{2}\eta\right\Vert _{H^{s}%
}^{2}+\\
+  &  \left\Vert V\right\Vert _{H^{s}}^{2}+\varepsilon(d-a)\left\Vert \nabla
V\right\Vert _{H^{s}}^{2}+\varepsilon^{2}(-a)d\left\Vert \nabla^{2}%
V\right\Vert _{H^{s}}^{2},
\end{align*}
where $\nabla^{2}\eta=\left(  \partial_{ij}^{2}\eta\right)  _{i,j}$ and
$\nabla^{2}V=\left(  \partial_{ij}^{2}V^{k}\right)  _{i,j,k}$. Clearly,
$H^{s_{1}}\left(  \mathbb{R}^{n}\right)  \times\left(  H^{s_{2}}\left(
\mathbb{R}^{n}\right)  \right)  ^{n}$ endowed with $\left\Vert \left(
\cdot,\cdot\right)  \right\Vert _{s}$ is a Banach space which is continuously
imbedded in $L^{2}\left(  \mathbb{R}^{n}\right)  \times\left(  L^{2}\left(
\mathbb{R}^{n}\right)  \right)  ^{n}$.

Our approach is based on an energy method applied to spectrally localized
equations. We first derive a priori estimates and we establish local existence
and uniqueness of solutions for the general $abcd$ system. Before we state the
main result let us formalize the notion of long time existence of solutions
for \eqref{eq2} in the next definitions.

\begin{definition}
Let $T>0$ a positive real number. We will say that $T$ is bounded from below
by a $O\left(  \frac{1}{\varepsilon}\right)  $-order quantity if there exists
another positive real number $C$, independent of $\varepsilon$ such that:
\[
T\geq\frac{C}{\varepsilon}.
\]

\end{definition}

Let us consider a Banach space $\left(  X,\left\Vert \cdot\right\Vert
_{X}\right)  $ which is continuously imbedded in $L^{2}\left(  \mathbb{R}%
^{n}\right)  \times\left(  L^{2}\left(  \mathbb{R}^{n}\right)  \right)  ^{n}$.

\begin{definition}
\label{def}Let us consider $W\in\left(  H^{1}\left(  \mathbb{R}^{n}\right)
\right)  ^{n}$. We say that we can establish long time existence and
uniqueness of solutions for the equation \eqref{eq2} in $X$ if for any
$\left(  \eta_{0},V_{0}\right)  \in X$ there exists a positive time $T>0$, an
unique solution\footnote{The solution should be understood in the tempered
distribution sense. See Definition \ref{Definitie}.} $\left(  \eta,V\right)
\in\mathcal{C}\left(  \left[  0,T\right]  ,X\right)  $ of \eqref{eq2} and a
function $F:\left(  0,+\infty\right)  \mathbb{\rightarrow}\left(
0,+\infty\right)  $ independent of $\varepsilon$ such that:%
\[
T\geq\frac{F\left(  \left\Vert \left(  \eta_{0},V_{0}\right)  \right\Vert
_{X}\right)  }{\varepsilon}.
\]

\end{definition}

\begin{remark}
Of course, the function $F$ appearing in the preceding definition is allowed
to depend on $a$,$b$,$c$,$d$,$W$ and on the dimension $n$.
\end{remark}

We are now in the position of stating our long time existence result:

\begin{theorem}
\label{Teorema2}Let $a,b,c,d$ as in \eqref{1111} excluding the two
cases\ \eqref{abcd}, $b+d>0$. Let us consider $s$ such that $s>\frac{n}{2}+1$
with $n\geq1$. Let us also consider $s_{1}$, $s_{2}$ and $s_{3}$ defined by
\eqref{relatie} and $W\in\left(  H^{s_{3}}\right)  ^{n}$. Then, we can
establish long time existence and uniqueness of solutions for the equation
\eqref{eq2} in $H^{s_{1}}\times\left(  H^{s_{2}}\right)  ^{n}$. Moreover, if
we denote by $T\left(  \eta_{0},V_{0}\right)  $, the maximal time of existence
then there exists some $T\in\lbrack0,T\left(  \eta_{0},V_{0}\right)  )$ which
is bounded from below by an $O\left(  \frac{1}{\varepsilon}\right)  $-order
quantity and a function $G:\mathbb{R\rightarrow R}$ such that for all
$t\in\left[  0,T\right]  $ we have:
\[
\left\Vert \left(  \eta,V\right)  \right\Vert _{s}\leq G\left(  \left\Vert
\left(  \eta_{0},V_{0}\right)  \right\Vert _{s}\right)  ,
\]
where $G$ may depend on $a,b,c,d,s,n$ but not on $\varepsilon$.
\end{theorem}

Theorem \ref{Teorema2} is the consequence of a more general result that we
obtain later in this paper. In fact, our method enables us (without any extra
effort) to establish long time existence and uniqueness of solutions in the
more general Besov spaces, thus achieving the critical regularity $s=\frac
{n}{2}+1$, see Theorem \ref{Teo}.

The rest of the paper is organized as follows. In Section \ref{S2} we
establish all the basic energy estimates that we will need in order to prove
Theorem \ref{Teorema2}. In Section \ref{S3} we prove that \eqref{eq2} admits
an unique solution and we establish an explosion criterion. The method used to
construct the solution assures an existence time that is bounded below by a
quantity of order $O\left(  \frac{1}{\sqrt{\varepsilon}}\right)  $. Finally in
Section \ref{iar} we prove Theorem \ref{Teorema2} namely we show that the
solution of \eqref{eq2} constructed in Section \ref{S3} persists on a time of
order $O\left(  \frac{1}{\varepsilon}\right)  $. The proof will be a
by-product of some refined energy estimates that we prove in Section \ref{S2}
and the explosion criteria established in Section \ref{S3}. In Section
\ref{sub1} we discuss about the cases \eqref{abcd}. We end the paper with
Section \ref{sub2} where we discuss the possibility of applying our method to
the $abcd$ systems in the case of a general bottom topography derived in
\cite{Chaz}, Section $2$.

\subsection{Notations}

Because our proof makes use of elementary tensor calculus let us introduce
some basic notations. For any vector field $U:\mathbb{R}^{n}\rightarrow
\mathbb{R}^{n}$ we denote by $\nabla U:\mathbb{R}^{n}\rightarrow
\mathcal{M}_{n}(\mathbb{R})$ and by $\nabla^{t}U:\mathbb{R}^{n}\rightarrow
\mathcal{M}_{n}(\mathbb{R})$ the $n\times n$ matrices defined by:%
\begin{align*}
\left(  \nabla U\right)  _{ij}  &  =\partial_{i}U^{j},\\
\left(  \nabla^{t}U\right)  _{ij}  &  =\partial_{j}U^{i}.
\end{align*}
In the same manner we define $\nabla^{2}U:\mathbb{R}^{n}\rightarrow
\mathbb{R}^{n}\times\mathbb{R}^{n}\times\mathbb{R}^{n}$ as:%
\[
\left(  \nabla^{2}U\right)  _{ijk}=\partial_{ij}^{2}U^{k}.
\]
We will suppose that all vectors appearing are column vectors and thus the
(classical) product between a matrix field $A$ and a vector field $U$ will be
the vector\footnote{From now on we will use the Einstein summation convention
over repeted indices.}:%
\[
\left(  AU\right)  ^{i}=A_{ij}U^{j}.
\]
We will often write the contraction operation between $\nabla^{2}U$ and a
vector field $V$ by%
\[
\left(  \nabla^{2}U:V\right)  _{ij}=\partial_{ij}^{2}U^{k}V^{k}%
\]
If $U,V:\mathbb{R}^{n}\rightarrow\mathbb{R}^{n}$ are two vector fields and
$A,B:\mathbb{R}^{n}\rightarrow\mathcal{M}_{n}(\mathbb{R})$ two matrix fields
we denote:%
\begin{align*}
UV  &  =U^{i}V^{i},\text{ }A:B=A_{ij}B_{ij},\\
\left\langle U,V\right\rangle _{L^{2}}  &  =\int U^{i}V^{i},\text{
}\left\langle A,B\right\rangle _{L^{2}}=\int A_{ij}B_{ij}\\
\left\Vert U\right\Vert _{L^{2}}^{2}  &  =\left\langle U,U\right\rangle
_{L^{2}},\text{ }\left\Vert A\right\Vert _{L^{2}}^{2}=\left\langle
A,A\right\rangle _{L^{2}}\\
\left\Vert \nabla^{2}U\right\Vert _{L^{2}}^{2}  &  =\int\nabla U:\nabla
U=\int\left(  \partial_{ij}U^{k}\right)  ^{2}%
\end{align*}
Also, the tensorial product of two vector fields $U,V$ is defined as the
matrix field $U\otimes V:\mathbb{R}^{n}\rightarrow\mathcal{M}_{n}(\mathbb{R})$
given by:%
\[
\left(  U\otimes V\right)  _{ij}=U^{i}V^{j}.
\]
We have the following derivation rule: if $u$ is a scalar field, $U$, $V$ are
vector fields and $A:\mathbb{R}^{n}\rightarrow\mathcal{M}_{n}(\mathbb{R})$
then:%
\begin{align*}
\nabla\operatorname{div}\left(  uU\right)   &  =\nabla^{2}uU+\nabla
u\operatorname{div}U+\nabla U\nabla u+u\nabla\operatorname{div}U\\
\nabla\left(  UV\right)   &  =\nabla UV+\nabla VU\\
\nabla(AV)  &  =\nabla^{2}A:V+\nabla VA^{t}%
\end{align*}
If we suppose that $\operatorname{curl}V=0$ then the following integration by
parts identity holds true:%
\begin{align}
\left\langle \nabla V:U,\nabla V\right\rangle _{L^{2}}  &  =\int\left(  \nabla
V:U\right)  :\nabla V=\int\partial_{ij}^{2}V^{k}U^{k}\partial_{i}V^{j}%
=\int\partial_{ik}^{2}V^{j}U^{k}\partial_{i}V^{j}\nonumber\\
&  =-\frac{1}{2}\int\partial_{k}U^{k}\left(  \partial_{i}V^{j}\right)
^{2}=-\frac{1}{2}\int\operatorname{div}U(\nabla V:\nabla V). \label{identity}%
\end{align}

\bigskip Let $\mathcal{C}$ be the annulus $\{\xi\in\mathbb{R}^{n}%
:3/4\leq\left\vert \xi\right\vert \leq8/3\}$. Let us choose two radial
functions $\chi\in\mathcal{D}(B(0,4/3))$ and $\varphi\in\mathcal{D(C)}$ valued
in the interval $\left[  0,1\right]  $ and such that:%
\[
\forall\xi\in\mathbb{R}^{n}\text{, \ }\chi(\xi)+\sum_{j\geq0}\varphi(2^{-j}%
\xi)=1\text{.}%
\]
Let us denote by $h=\mathcal{F}^{-1}\varphi$ and $\tilde{h}=\mathcal{F}%
^{-1}\chi$. For all $u\in\mathcal{S}^{\prime}$, the nonhomogeneous dyadic
blocks are defined as follows:%
\begin{gather}
\Delta_{j}u=0\text{ \ if \ }j\leq-2,\nonumber\\
\Delta_{-1}u=\chi\left(  D\right)  u=\tilde{h}\star u,\label{diadice}\\
\Delta_{j}u=\varphi\left(  2^{-j}D\right)  u=2^{jd}\int_{\mathbb{R}^{n}%
}h\left(  2^{q}y\right)  u\left(  x-y\right)  dy\text{ \ if \ }j\geq
0.\nonumber
\end{gather}
The high frequency cut-off operator $S_{j}$ is defined by%
\[
S_{j}u=\sum_{k\leq j-1}\Delta_{k}u.
\]
Let us define now the nonhomogeneous Besov spaces.

\begin{definition}
Let $s\in\mathbb{R}$, $r\in\left[  1,\infty\right]  $. The Besov space
$B_{2,r}^{s}$ is the set of tempered distributions $u\in\mathcal{S}^{\prime}$
such that:%
\[
\left\Vert u\right\Vert _{B_{2,r}^{s}}:=\left\Vert \left(  2^{js}\left\Vert
\Delta_{j}u\right\Vert _{L^{2}}\right)  _{j\in\mathbb{Z}}\right\Vert
_{\ell^{r}(\mathbb{Z)}}<\infty.
\]

\end{definition}

Let us mention that $H^{s}\left(  \mathbb{R}^{n}\right)  =B_{2,2}^{s}\left(
\mathbb{R}^{n}\right)  $ and that we have the following continuous\ embedding%
\[
B_{2,1}^{s}\hookrightarrow H^{s}\hookrightarrow B_{2,\infty}^{s}%
\hookrightarrow H^{s^{\prime}},
\]
for all $s^{\prime}<s$. Some basic properties about Besov spaces can be found
in the Appendix. For more details and full proofs we refer to \cite{Dan1}.
\bigskip

Let us consider $\varepsilon\leq1$ and $s>0$, $r\in\left[  1,\infty\right]  $.
For all $\left(  \eta,V\right)  \in B_{2,r}^{s_{1}}\left(  \mathbb{R}%
^{n}\right)  \times\left(  B_{2,r}^{s_{2}}\left(  \mathbb{R}^{n}\right)
\right)  ^{n}$ we introduce the following quantities:%
\begin{align}
U_{j}^{2}  &  =U_{j}^{2}\left(  \eta,V\right)  =\int\eta_{j}^{2}%
+\varepsilon\left(  b-c\right)  \left\vert \nabla\eta_{j}\right\vert
^{2}+\varepsilon^{2}(-c)b\left(  \nabla^{2}\eta_{j}:\nabla^{2}\eta_{j}\right)
\label{Uj}\\
&  +\int V_{j}^{2}+\varepsilon\left(  d-a\right)  \left(  \nabla V_{j}:\nabla
V_{j}\right)  +\varepsilon^{2}(-a)d\left(  \nabla^{2}V_{j}:\nabla^{2}%
V_{j}\right)  ,\nonumber
\end{align}
and%

\begin{align}
U_{s}^{2}  &  =U_{s}^{2}\left(  \eta,V\right)  =\left\Vert \eta\right\Vert
_{B_{2,r}^{s}}^{2}+\varepsilon\left(  b-c\right)  \left\Vert \nabla
\eta\right\Vert _{B_{2,r}^{s}}^{2}+\varepsilon^{2}\left(  -c\right)
b\left\Vert \nabla^{2}\eta\right\Vert _{B_{2,r}^{s}}^{2}\label{Us}\\
&  +\left\Vert V\right\Vert _{B_{2,r}^{s}}^{2}+\varepsilon\left(  d-a\right)
\left\Vert \nabla V\right\Vert _{B_{2,r}^{s}}^{2}+\varepsilon^{2}\left(
-a\right)  d\left\Vert \nabla^{2}V\right\Vert _{B_{2,r}^{s}}^{2}.\nonumber
\end{align}
where $\left(  \eta_{j},V_{j}\right)  :=(\Delta_{j}\eta,\Delta_{j}V)$ are the
frequency-localized dyadic blocks defined by relation \eqref{diadice}. It is
easy to check that $U_{s}\left(  \eta,V\right)  $ is a norm on the space
$B_{2,r}^{s_{1}}\left(  \mathbb{R}^{n}\right)  \times\left(  B_{2,r}^{s_{2}%
}\left(  \mathbb{R}^{n}\right)  \right)  ^{n}$. Also, it transpires that:%
\[
\left\Vert \left(  2^{js}U_{j}\right)  _{j\in\mathbb{Z}}\right\Vert _{\ell
^{r}}\approx U_{s}.
\]
Some times we will also use the notation
\[
\left\Vert \left(  \eta,V\right)  \right\Vert _{s}=U_{s}\left(  \eta,V\right)
,
\]
for all $\left(  \eta,V\right)  \in B_{2,r}^{s_{1}}\left(  \mathbb{R}%
^{n}\right)  \times\left(  B_{2,r}^{s_{2}}\left(  \mathbb{R}^{n}\right)
\right)  ^{n}$ . We observe that $\left(  B_{2,r}^{s_{1}}\left(
\mathbb{R}^{n}\right)  \times\left(  B_{2,r}^{s_{2}}\left(  \mathbb{R}%
^{n}\right)  \right)  ^{n},\left\Vert \left(  \cdot,\cdot\right)  \right\Vert
_{s}\right)  $ is a Banach space. In order to ease the notations we will
rather write $B_{2,r}^{s}$ instead of $B_{2,r}^{s}\left(  \mathbb{R}%
^{n}\right)  $. Another quantity that will play an important role in the
following analysis is:%
\begin{gather}
N_{j}^{2}=N_{j}^{2}\left(  \eta,V\right)  =\int\left(  1+\varepsilon\left\Vert
\eta\right\Vert _{L^{\infty}}\right)  \eta_{j}^{2}+\varepsilon\left(
b-c\right)  \left(  1+\varepsilon\left\Vert \eta\right\Vert _{L^{\infty}%
}\right)  \left\vert \nabla\eta_{j}\right\vert ^{2}\nonumber\\
+\int\varepsilon^{2}(-c)b\left(  1+\varepsilon\left\Vert \eta\right\Vert
_{L^{\infty}}\right)  \left(  \nabla^{2}\eta_{j}:\nabla^{2}\eta_{j}\right)
\nonumber\\
+\int\left(  1+\varepsilon\eta+\varepsilon\left\Vert \eta\right\Vert
_{L^{\infty}}\right)  V_{j}^{2}+\varepsilon\left(  d-a+d\varepsilon
\eta+d\varepsilon\left\Vert \eta\right\Vert _{L^{\infty}}\right)  \left(
\nabla V_{j}:\nabla V_{j}\right) \nonumber\\
+\int\varepsilon^{2}(-a)d\left(  1+\varepsilon\left\Vert \eta\right\Vert
_{L^{\infty}}\right)  \left(  \nabla^{2}V_{j}:\nabla^{2}V_{j}\right)  ,
\label{Nj}%
\end{gather}
which satisfies:%
\begin{equation}
U_{j}\left(  \eta,V\right)  \lesssim N_{j}\left(  \eta,V\right)
\lesssim\left(  1+2\varepsilon\left\Vert \eta\right\Vert _{L^{\infty}}\right)
^{\frac{1}{2}}U_{j}\left(  \eta,V\right)  \label{relNj}%
\end{equation}
Denoting by
\begin{equation}
N_{s}=N_{s}\left(  \eta,V\right)  =\left\Vert \left(  2^{js}N_{j}\left(
\eta,V\right)  \right)  _{j\in\mathbb{Z}}\right\Vert _{\ell^{r}\left(
\mathbb{Z}\right)  } \label{Ns}%
\end{equation}
we obviously have that%
\begin{equation}
U_{s}\left(  \eta,V\right)  \lesssim N_{s}\left(  \eta,V\right)
\lesssim\left(  1+2\varepsilon\left\Vert \eta\right\Vert _{L^{\infty}}\right)
^{\frac{1}{2}}U_{s}\left(  \eta,V\right)  . \label{relNs}%
\end{equation}

\section{Energy-type identities\label{S2}}

We begin by localizing equation \eqref{eq2} in Fourier space thus obtaining
that:%
\begin{equation}
\left\{
\begin{array}
[c]{l}%
\left(  I-\varepsilon b\Delta\right)  \partial_{t}\eta_{j}+\operatorname{div}%
V_{j}+a\varepsilon\operatorname{div}\Delta V_{j}+\varepsilon W\nabla\eta
_{j}+\varepsilon V\nabla\eta_{j}+\varepsilon\eta\operatorname{div}%
V_{j}=\varepsilon R_{1j}\\
\left(  I-\varepsilon d\Delta\right)  \partial_{t}V_{j}+\nabla\eta
_{j}+c\varepsilon\nabla\Delta\eta_{j}+\varepsilon\nabla V_{j}W+\varepsilon
\nabla V_{j}V=\varepsilon R_{2j}\\
\eta_{j|t=0}=\Delta_{j}\eta_{0}\text{ , }V_{j|t=0}=\Delta_{j}V_{0}%
\end{array}
\right.  \label{eq3}%
\end{equation}
where the remainder terms are given by\footnote{From now on, we agree that if
$A$, $B$ are two operator then the operator $\left[  A,B\right]  $ is given
by:%
\[
\left[  A,B\right]  =AB-BA.
\]
}:%
\begin{equation}
\left\{
\begin{array}
[c]{l}%
R_{1j}=\left[  W,\Delta_{j}\right]  \nabla\eta+\left[  V,\Delta_{j}\right]
\nabla\eta+\left[  \eta,\Delta_{j}\right]  \operatorname{div}V\\
R_{2j}=\left[  W,\Delta_{j}\right]  \nabla V+\left[  V,\Delta_{j}\right]
\nabla V-\frac{1}{2}\nabla\Delta_{j}\left(  \left\vert W\right\vert
^{2}\right)  -\Delta_{j}\left(  \nabla WV\right)  .
\end{array}
\right.  \label{eq4}%
\end{equation}
Let us establish our first useful identity. We multiply the first equation in
\eqref{eq3} by $\eta_{j}$ and the second one with $V_{j}$ and by adding them
up and integrating in space we get that\footnote{Observe that here we use the
fact that $\operatorname{curl}V=0$. Indeed, under \eqref{curl} we have%
\[
\int\nabla V_{j}VV_{j}=\frac{1}{2}\int V\nabla\left\vert V_{j}\right\vert
^{2}=-\frac{1}{2}\int\operatorname{div}V\left\vert V_{j}\right\vert ^{2}.
\]
}:%
\begin{gather}
\frac{1}{2}\partial_{t}\left(  \left\Vert \eta_{j}\right\Vert _{L^{2}}%
^{2}+\varepsilon b\left\Vert \nabla\eta_{j}\right\Vert _{L^{2}}^{2}+\left\Vert
V_{j}\right\Vert _{L^{2}}^{2}+\varepsilon d\left\Vert \nabla V_{j}\right\Vert
_{L^{2}}^{2}\right)  +a\varepsilon\int\eta_{j}\operatorname{div}\Delta
V_{j}+c\varepsilon\int V_{j}\nabla\Delta\eta_{j}\label{eq4'}\\
+\varepsilon\int\eta\eta_{j}\operatorname{div}V_{j}=\frac{\varepsilon}{2}%
\int\operatorname{div}V\left(  \eta_{j}^{2}+V_{j}^{2}\right)  +\varepsilon\int
R_{1j}\eta_{j}+\varepsilon\int R_{2j}V_{j}.\nonumber
\end{gather}
Let us denote by $T_{1}$ the right hand side of the above identity:%
\begin{equation}
T_{1}=\frac{\varepsilon}{2}\int\operatorname{div}V\left(  \eta_{j}^{2}%
+V_{j}^{2}\right)  +\varepsilon\int R_{1j}\eta_{j}+\varepsilon\int R_{2j}%
V_{j}. \label{T1}%
\end{equation}
Next, we wish to derive similar identities involving the quantities
$\nabla\eta_{j},\nabla^{2}\eta_{j}$ and $\nabla V_{j},\nabla^{2}V_{j}$. In
order to do so, let us observe that applying $\nabla$ to the first equation in
\eqref{eq2} gives us:%
\[
\left(  I-\varepsilon b\Delta\right)  \partial_{t}\nabla\eta+\nabla
\operatorname{div}V+a\varepsilon\nabla\operatorname{div}\Delta V+\varepsilon
\nabla W\nabla\eta+\varepsilon\nabla^{2}\eta W+\varepsilon\nabla
\operatorname{div}\left(  \eta V\right)  =0
\]
and that by applying $\Delta_{j}$ we end up with:%
\begin{gather}
\left(  I-\varepsilon b\Delta\right)  \partial_{t}\nabla\eta_{j}%
+\nabla\operatorname{div}V_{j}+a\varepsilon\nabla\operatorname{div}\Delta
V_{j}+\varepsilon\nabla^{2}\eta_{j}W+\varepsilon\nabla^{2}\eta_{j}%
V+\varepsilon\eta\nabla\operatorname{div}V_{j}\label{eq5''}\\
+\varepsilon\operatorname{div}V_{j}\nabla\eta+\varepsilon\nabla V_{j}%
\nabla\eta=-\varepsilon\Delta_{j}\left(  \nabla W\nabla\eta\right)
+\varepsilon R_{3j},\nonumber
\end{gather}
where%
\[
R_{3j}=\left[  W,\Delta_{j}\right]  \nabla^{2}\eta+\left[  V,\Delta
_{j}\right]  \nabla^{2}\eta+\left[  \eta,\Delta_{j}\right]  \nabla
\operatorname{div}V+[\nabla\eta,\Delta_{j}]\operatorname{div}V+\left[
\nabla\eta,\Delta_{j}\right]  \nabla V.
\]
We multiply \eqref{eq5''} with $-c\varepsilon\nabla\eta_{j}$ and by
integration we get\footnote{Here, we use the fact that $\operatorname{div}%
W=0$. We get that:%
\[
\varepsilon\int W\nabla\eta_{j}\eta_{j}=-\frac{\varepsilon}{2}\int
\operatorname{div}W\eta_{j}^{2}=0.
\]
}:%
\begin{gather}
\frac{-c\varepsilon}{2}\partial_{t}\left(  \int\left(  \left\vert \nabla
\eta_{j}\right\vert ^{2}+\varepsilon b\nabla^{2}\eta_{j}:\nabla^{2}\eta
_{j}\right)  \right)  -c\varepsilon\int\nabla\operatorname{div}V_{j}\nabla
\eta_{j}-ac\varepsilon^{2}\int\nabla\operatorname{div}\Delta V_{j}\nabla
\eta_{j}\label{eq5'}\\
-c\varepsilon^{2}\int\eta\nabla\operatorname{div}V_{j}\nabla\eta_{j}%
=T_{2}\nonumber
\end{gather}
with $T_{2}$ given by%
\begin{gather}
T_{2}=c\varepsilon^{2}\int(\nabla^{2}\eta_{j}W)\nabla\eta_{j}+c\varepsilon
^{2}\int(\nabla^{2}\eta_{j}V)\nabla\eta_{j}\nonumber\\
+c\varepsilon^{2}\int\operatorname{div}V_{j}\nabla\eta\nabla\eta
_{j}+c\varepsilon^{2}\int(\nabla V_{j}\nabla\eta)\nabla\eta_{j}+c\varepsilon
^{2}\int\Delta_{j}\left(  \nabla W\nabla\eta\right)  \nabla\eta_{j}%
-c\varepsilon^{2}\int R_{3j}\nabla\eta_{j}. \label{T2}%
\end{gather}
We proceed similarly with the second equation in \eqref{eq2} and we obtain:%
\begin{gather*}
\left(  I-\varepsilon d\Delta\right)  \partial_{t}\nabla V+\nabla^{2}%
\eta+\varepsilon c\nabla^{2}\Delta\eta+\varepsilon\frac{1}{2}\nabla
^{2}\left\vert W\right\vert ^{2}+\varepsilon\left(  \nabla^{2}W:V+\nabla
V\nabla^{t}W\right) \\
+\varepsilon\left(  \nabla^{2}V:W+\nabla W\nabla^{t}V\right)  +\varepsilon
\left(  \nabla^{2}V:V+\nabla V\nabla^{t}V\right)  =0.
\end{gather*}
We localize the last equation and we get that:%
\begin{gather}
\left(  I-\varepsilon d\Delta\right)  \partial_{t}\nabla V_{j}+\nabla^{2}%
\eta_{j}+\varepsilon c\nabla^{2}\Delta\eta_{j}+\varepsilon\frac{1}{2}%
\Delta_{j}\nabla^{2}\left\vert W\right\vert ^{2}+\varepsilon\Delta_{j}\left(
\nabla^{2}W:V\right) \label{eq6}\\
+\varepsilon\Delta_{j}\left(  \nabla V\nabla^{t}W\right)  +\varepsilon
\nabla^{2}V_{j}:W+\varepsilon\Delta_{j}\left(  \nabla W\nabla^{t}V\right)
+\varepsilon\nabla^{2}V_{j}:V+\varepsilon\nabla V_{j}\nabla^{t}V=\varepsilon
R_{4j},\nonumber
\end{gather}
where:%
\[
R_{4j}=\left[  W,\Delta_{j}\right]  \nabla^{2}V+\left[  V,\Delta_{j}\right]
\nabla^{2}V+[\nabla^{t}V,\Delta_{j}]\nabla V.
\]
We contract \eqref{eq6} with $-a\varepsilon\nabla V_{j}$ and by integration we
get that:%
\begin{equation}
\frac{-a\varepsilon}{2}\partial_{t}\left(  \int\nabla V_{j}:\nabla
V_{j}+\varepsilon d\nabla^{2}V_{j}:\nabla^{2}V_{j}\right)  -a\varepsilon
\int\nabla^{2}\eta_{j}:\nabla V_{j}-ac\varepsilon^{2}\int\nabla^{2}\Delta
\eta_{j}:\nabla V_{j}=T_{3} \label{eq6''}%
\end{equation}
with $T_{3}$ given by:%
\begin{gather}
T_{3}=a\varepsilon^{2}\int\left(  \nabla^{2}V_{j}:W\right)  :\nabla
V_{j}+a\varepsilon^{2}\int\left(  \nabla^{2}V_{j}:V\right)  :\nabla
V_{j}+a\varepsilon^{2}\int\nabla V_{j}\nabla^{t}V:\nabla V_{j}-a\varepsilon
^{2}\int R_{4j}:\nabla V_{j}\nonumber\\
+a\varepsilon^{2}\int\left(  \frac{1}{2}\Delta_{j}\nabla^{2}\left\vert
W\right\vert ^{2}+\Delta_{j}\left(  \nabla W\nabla^{t}V\right)  +\Delta
_{j}\left(  \nabla^{2}W:V\right)  +\Delta_{j}\left(  \nabla V\nabla
^{t}W\right)  \right)  :\nabla V_{j} \label{T3}%
\end{gather}
Let us add up identities \eqref{eq5'} and \eqref{eq6''} to get
that\footnote{Observe that by integration by parts we get
\[
-ac\varepsilon^{2}\int\nabla\operatorname{div}\Delta V_{j}:\nabla\eta
_{j}-ac\varepsilon^{2}\int\nabla^{2}\Delta\eta_{j}:\nabla V_{j}=0.
\]
}:%
\begin{gather*}
\partial_{t}\left(  \int\left\{  \frac{-c\varepsilon}{2}\left(  \left\vert
\nabla\eta_{j}\right\vert ^{2}+\varepsilon b\nabla^{2}\eta_{j}:\nabla^{2}%
\eta_{j}\right)  +\frac{-a\varepsilon}{2}\left(  \nabla V_{j}:\nabla
V_{j}+\varepsilon d\nabla^{2}V_{j}:\nabla^{2}V_{j}\right)  \right\}  \right)
\\
-a\varepsilon\int\nabla^{2}\eta_{j}:\nabla V_{j}-c\varepsilon\int
\nabla\operatorname{div}V_{j}\nabla\eta_{j}-c\varepsilon^{2}\int\eta
\nabla\operatorname{div}V_{j}\nabla\eta_{j}=T_{2}+T_{3}.
\end{gather*}
Finally add up \eqref{eq4'} to the last identity in order to
obtain\footnote{Observe that the terms $a\varepsilon\int\eta_{j}%
\operatorname{div}\Delta V_{j}-a\varepsilon\int\nabla^{2}\eta_{j}:\nabla
V_{j}$ \ and $c\varepsilon\int V_{j}\nabla\Delta\eta_{j}-c\varepsilon
\int\nabla\operatorname{div}V_{j}:\nabla\eta_{j}$ are both $0$.}:%
\begin{gather}
\frac{1}{2}\partial_{t}\left(  \int\eta_{j}^{2}+\varepsilon\left(  b-c\right)
\left\vert \nabla\eta_{j}\right\vert ^{2}+\varepsilon^{2}(-c)b\left(
\nabla^{2}\eta_{j}:\nabla^{2}\eta_{j}\right)  \right) \nonumber\\
+\frac{1}{2}\partial_{t}\left(  \int\left\vert V_{j}\right\vert ^{2}%
+\varepsilon\left(  d-a\right)  \left(  \nabla V_{j}:\nabla V_{j}\right)
+\varepsilon^{2}(-a)d\left(  \nabla^{2}V_{j}:\nabla^{2}V_{j}\right)  \right)
\nonumber\\
+\varepsilon\int\eta\eta_{j}\operatorname{div}V_{j}-c\varepsilon^{2}\int
\eta\nabla\operatorname{div}V_{j}\nabla\eta_{j}=T_{1}+T_{2}+T_{3}. \label{eq7}%
\end{gather}

\subsection{Refined energy-type identities\label{refined}}

As we have already seen in \eqref{Hamil}, the $abcd$ system possesses a
formally conserved energy. The types of estimates that we establish in this
section, resemble very much to this conserved energy and will be the
equivalent of the ones obtained in \cite{Bona3} pages $617$ and $625$, with
$\Delta_{j}\eta$ and $\Delta_{j}V$ instead of $\eta$ and $V$. Of course this
is the key ingredient for obtaining long time existence results that allow
initial data to lie in larger spaces.

Having established identity \eqref{eq7} we observe that the terms
$\varepsilon\int\eta\eta_{j}\operatorname{div}V_{j}$ and $-c\varepsilon
^{2}\int\eta\nabla\operatorname{div}V_{j}\nabla\eta_{j}$ prevent us from
directly applying a Gronwall type argument and establishing long time
existence. In order to repair this inconvenience let us multiply the second
equation of \eqref{eq3} with $\varepsilon\eta V_{j}$. We thus get:%
\begin{align*}
&  \left\langle \varepsilon\left(  I-d\Delta\right)  \partial_{t}V_{j},\eta
V_{j}\right\rangle _{L^{2}}+\varepsilon\int\eta\nabla\eta_{j}V_{j}%
+c\varepsilon^{2}\int\eta\nabla\Delta\eta_{j}V_{j}\\
&  =-\varepsilon^{2}\int\left(  \nabla V_{j}W\right)  (\eta V_{j}%
)-\varepsilon^{2}\int\left(  \nabla V_{j}V\right)  \left(  \eta V_{j}\right)
+\varepsilon^{2}\int\eta R_{2j}V_{j}.
\end{align*}
Observing that the first term writes:%
\begin{align*}
\left\langle \varepsilon\left(  I-d\Delta\right)  \partial_{t}V_{j},\eta
V_{j}\right\rangle _{L^{2}}  &  =\frac{1}{2}\partial_{t}\left(  \int
\varepsilon\eta\left\vert V_{j}\right\vert ^{2}+\varepsilon^{2}d\eta\nabla
V_{j}:\nabla V_{j}\right)  -\frac{\varepsilon}{2}\int\partial_{t}%
\eta\left\vert V_{j}\right\vert ^{2}\\
-\varepsilon^{2}\frac{d}{2}\int\partial_{t}\eta\nabla V_{j}  &  :\nabla
V_{j}+\varepsilon^{2}d\left\langle \partial_{t}\nabla V_{j},\nabla\eta\otimes
V_{j}\right\rangle _{L^{2}}%
\end{align*}
we write that
\begin{equation}
\frac{1}{2}\partial_{t}\left(  \int\varepsilon\eta\left\vert V_{j}\right\vert
^{2}+\varepsilon^{2}d\eta\nabla V_{j}:\nabla V_{j}\right)  +\varepsilon
\int\eta\nabla\eta_{j}V_{j}+c\varepsilon^{2}\int\eta\nabla\Delta\eta_{j}%
V_{j}=T_{4} \label{eq8}%
\end{equation}
with $T_{4}$ given by:%
\begin{gather}
T_{4}=-\varepsilon^{2}\int\left(  \nabla V_{j}W\right)  \eta V_{j}%
-\varepsilon^{2}\int\left(  \nabla V_{j}V\right)  \eta V_{j}+\varepsilon
^{2}\int\eta R_{2j}V_{j}+\nonumber\\
+\frac{\varepsilon}{2}\int\partial_{t}\eta\left\vert V_{j}\right\vert
^{2}+\varepsilon^{2}\frac{d}{2}\int\partial_{t}\eta\nabla V_{j}:\nabla
V_{j}-\varepsilon^{2}d\left\langle \partial_{t}\nabla V_{j},\nabla\eta\otimes
V_{j}\right\rangle _{L^{2}}. \label{T4}%
\end{gather}
We add up \eqref{eq7} and \eqref{eq8} in order to obtain:%
\begin{gather}
\frac{1}{2}\partial_{t}\left(  \int\eta_{j}^{2}+\varepsilon\left(  b-c\right)
\left\vert \nabla\eta_{j}\right\vert ^{2}+\varepsilon^{2}(-c)b\left(
\nabla^{2}\eta_{j}:\nabla^{2}\eta_{j}\right)  \right) \nonumber\\
+\frac{1}{2}\partial_{t}\left(  \int\left(  1+\varepsilon\eta\right)
\left\vert V_{j}\right\vert ^{2}+\varepsilon\left(  d-a+d\varepsilon
\eta\right)  \left(  \nabla V_{j}:\nabla V_{j}\right)  +\varepsilon
^{2}(-a)d\left(  \nabla^{2}V_{j}:\nabla^{2}V_{j}\right)  \right) \nonumber\\
=T_{0}+T_{1}+T_{2}+T_{3}+T_{4}, \label{eq9}%
\end{gather}
with%
\[
T_{0}=\varepsilon\int\nabla\eta\eta_{j}V_{j}-c\varepsilon^{2}\int\nabla
\eta\Delta\eta_{j}V_{j}+c\varepsilon^{2}\int\nabla\eta\nabla\eta
_{j}\operatorname{div}V_{j}.
\]

\subsection{Estimates for the $T_{i}$'s\label{Sec}}

Having established the energy identity \eqref{eq9}, we proceed by conveniently
bounding the RHS term. We want to obtain a bound of the form:%
\[
T_{0}+T_{1}+T_{2}+T_{3}+T_{4}\leq\varepsilon P(U_{j},U_{s})
\]
with $P\left(  x,y\right)  $ some polynomial function with coefficients not
depending on $\varepsilon$. We suppose that
\[
s>\frac{n}{2}+1\text{ or }s=\frac{n}{2}+1\text{ and }r=1\text{.}
\]
Moreover, $C>0$ will denote a generic positive constant depending only on the
dimension $n$ and on $s$. All the estimates established here are valid for the
$a,b,c,d$ parameters satisfying \eqref{1111} with $b+d>0$. The case:
\begin{equation}
a=d=0,\text{ }b>0,\text{ }c<0 \label{abcd1}%
\end{equation}
needs more attention and we will investigate it in Section \ref{Section}, thus
for the moment we do all the computation supposing that case \eqref{abcd1}
does nor occur.

We claim that we can bound $T_{0}+T_{1}+T_{2}+T_{3}+T_{4}$ in such a way that
the following estimate holds true:%
\begin{gather}
\frac{1}{2}\partial_{t}\left(  \int\eta_{j}^{2}+\varepsilon\left(  b-c\right)
\left\vert \nabla\eta_{j}\right\vert ^{2}+\varepsilon^{2}(-c)b\left(
\nabla^{2}\eta_{j}:\nabla^{2}\eta_{j}\right)  \right) \nonumber\\
+\frac{1}{2}\partial_{t}\left(  \int\left(  1+\varepsilon\eta\right)
V_{j}^{2}+\varepsilon\left(  d-a+d\varepsilon\eta\right)  \left(  \nabla
V_{j}:\nabla V_{j}\right)  +\varepsilon^{2}(-a)d\left(  \nabla^{2}V_{j}%
:\nabla^{2}V_{j}\right)  \right) \nonumber\\
\leq\varepsilon CU_{j}\left\{  U_{j}\left(  U_{s}+HU_{s}+U_{s}^{2}\right)
+c_{j}(t)(1+U_{s})\left(  H^{2}+HU_{s}+U_{s}^{2}\right)  \right\}  ,
\label{eq14}%
\end{gather}
where $C>0$ depends only on $a$,$b$,$c$,$d$,$n$ and $s$ but not on
$\varepsilon$, $\left(  c_{j}\left(  t\right)  \right)  _{j}$ is a sequence
with $\left(  2^{js}c_{j}\left(  t\right)  \right)  _{j}\in\ell^{r}%
(\mathbb{Z)}$, having norm $1$ and finally $H$ is defined by:%
\[
H=\left\Vert W\right\Vert _{B_{2,r}^{s}}+\left\Vert \nabla W\right\Vert
_{B_{2,r}^{s}}-sgn\left(  a\right)  \sqrt{\varepsilon}\left\Vert \nabla
^{2}W\right\Vert _{B_{2,r}^{s}}.
\]
Let us detail this below. Regarding $T_{0}$, we begin by observing that:%
\[
\varepsilon\int\nabla\eta\eta_{j}V_{j}\leq\varepsilon\left\Vert \nabla
\eta\right\Vert _{L^{\infty}}\left\Vert \eta_{j}\right\Vert _{L^{2}}\left\Vert
V_{j}\right\Vert _{L^{2}}\leq\varepsilon U_{j}^{2}U_{s}.
\]
Next, if $a=d=0$ then $b>0$ and we get that\footnote{From now on, we adopt the
convention $\frac{0}{0}=0$.}:
\begin{align*}
&  -c\varepsilon^{2}\int\Delta\eta\nabla\eta_{j}V_{j}+c\varepsilon^{2}%
\int\nabla\eta\nabla\eta_{j}\operatorname{div}V_{j}\\
&  =-c\varepsilon^{2}\int\Delta\eta\nabla\eta_{j}V_{j}-c\varepsilon^{2}%
\int\nabla^{2}\eta\nabla\eta_{j}V_{j}-c\varepsilon^{2}\int\nabla^{2}\eta
_{j}\nabla\eta V_{j}\\
&  \leq-c\varepsilon^{2}\left(  \left\Vert \Delta\eta\right\Vert _{L^{\infty}%
}\left\Vert \nabla\eta_{j}\right\Vert _{L^{2}}\left\Vert V_{j}\right\Vert
_{L^{2}}+\left\Vert \nabla^{2}\eta\right\Vert _{L^{\infty}}\left\Vert
\nabla\eta_{j}\right\Vert _{L^{2}}\left\Vert V_{j}\right\Vert _{L^{2}}\right.
\\
&  \text{ \ \ \ \ \ \ \ \ \ \ \ \ \ }\left.  +\left\Vert \nabla\eta\right\Vert
_{L^{\infty}}\left\Vert \nabla^{2}\eta_{j}\right\Vert \left\Vert
V_{j}\right\Vert _{L^{2}}\right) \\
&  \leq\varepsilon C\max\left(  \frac{-c}{b-c},\sqrt{\frac{-c}{b}}\right)
U_{j}^{2}U_{s}.
\end{align*}
If at least one of $a,d$ is not $0$ then we write:%
\begin{align*}
&  -c\varepsilon^{2}\int\nabla\eta\Delta\eta_{j}V_{j}+c\varepsilon^{2}%
\int\nabla\eta\nabla\eta_{j}\operatorname{div}V_{j}\\
&  =c\varepsilon^{2}\int\left(  \nabla^{2}\eta V_{j}\right)  \nabla\eta
_{j}+c\varepsilon^{2}\int\left(  \nabla V_{j}\nabla\eta\right)  \nabla\eta
_{j}++c\varepsilon^{2}\int\nabla\eta\nabla\eta_{j}\operatorname{div}V_{j}\\
&  \leq-c\varepsilon^{2}\left(  \left\Vert \nabla^{2}\eta\right\Vert
_{L^{\infty}}\left\Vert V_{j}\right\Vert _{L^{2}}\left\Vert \nabla\eta
_{j}\right\Vert _{L^{2}}+\left\Vert \nabla\eta\right\Vert _{L^{\infty}%
}\left\Vert \nabla V_{j}\right\Vert _{L^{2}}\left\Vert \nabla\eta
_{j}\right\Vert _{L^{2}}\right. \\
\text{ }  &  \left.  \text{ \ \ \ \ \ \ \ \ \ \ \ \ \ }+\left\Vert \nabla
\eta\right\Vert _{L^{\infty}}\left\Vert \nabla\eta_{j}\right\Vert _{L^{2}%
}\left\Vert \operatorname{div}V_{j}\right\Vert _{L^{2}}\right) \\
&  \leq\varepsilon C\max\left(  \frac{-c}{b-c},\frac{-c}{d-a}\right)
U_{j}^{2}U_{s}.
\end{align*}
We choose
\begin{equation}
C_{abcd}^{1}=\left\{
\begin{array}
[c]{l}%
\max\left(  \frac{-c}{b-c},\sqrt{\frac{-c}{b}}\right)  \text{ if }d-a=0,\\
\max\left(  \frac{-c}{b-c},\frac{-c}{d-a}\right)  \text{ if }d-a>0\text{.}%
\end{array}
\right.  \label{1}%
\end{equation}
Consequently, we have:%
\begin{equation}
T_{0}\leq\varepsilon CC_{abcd}^{1}U_{j}^{2}U_{s}. \label{estim1}%
\end{equation}

Next, let us analyze $T_{1}$. According to Proposition \ref{comut} and
Proposition \ref{propA} we get that:%
\begin{align*}
\left\Vert R_{1j}\right\Vert _{L^{2}}  &  \leq Cc_{j}^{1}\left(  t\right)
\left(  \left\Vert \nabla W\right\Vert _{B_{2,r}^{s-1}}\left\Vert
\eta\right\Vert _{B_{2,r}^{s}}+\left\Vert \nabla V\right\Vert _{B_{2,r}^{s-1}%
}\left\Vert \eta\right\Vert _{B_{2,r}^{s}}+\left\Vert \nabla\eta\right\Vert
_{B_{2,r}^{s-1}}\left\Vert V\right\Vert _{B_{2,r}^{s}}\right)  ,\\
\left\Vert R_{2j}\right\Vert _{L^{2}}  &  \leq Cc_{j}^{2}\left(  t\right)
\left(  \left(  \left\Vert \nabla W\right\Vert _{B_{2,r}^{s-1}}+\left\Vert
\nabla V\right\Vert _{B_{2,r}^{s-1}}\right)  \left\Vert V\right\Vert
_{B_{2,r}^{s}}+\left(  \left\Vert W\right\Vert _{B_{2,r}^{s}}+\left\Vert
V\right\Vert _{B_{2,r}^{s}}\right)  \left\Vert \nabla W\right\Vert
_{B_{2,r}^{s}}\right)  ,
\end{align*}
with $\left(  2^{js}c_{j}^{i}\left(  t\right)  \right)  \in\ell^{r}%
(\mathbb{Z)}$ , $i=1,2$, with norm $1$. In order to avoid mentioning each
time, from now on, \ for all natural number, $i\in\mathbb{N}$, $\left(
c_{j}^{i}\left(  t\right)  \right)  _{j\in\mathbb{Z}}$ will be a sequence such
that $\left(  2^{js}c_{j}^{i}\left(  t\right)  \right)  _{j\in\mathbb{Z}}%
\in\ell^{r}(\mathbb{Z)}$ with norm $1$. Recall that $H$ stands for the
following quantity:%
\begin{equation}
H=\left\Vert W\right\Vert _{B_{2,r}^{s}}+\left\Vert \nabla W\right\Vert
_{B_{2,r}^{s}}-sgn\left(  a\right)  \sqrt{\varepsilon}\left\Vert \nabla
^{2}W\right\Vert _{B_{2,r}^{s}} \label{relH}%
\end{equation}
and rewrite the previous inequalities (of course the constant $C=C\left(
n,s\right)  $ changes whenever necessary) as:%
\begin{align}
\left\Vert R_{1j}\right\Vert _{L^{2}}  &  \leq Cc_{j}^{1}\left(  t\right)
U_{s}\left(  H+U_{s}\right)  ,\label{R1}\\
\left\Vert R_{2j}\right\Vert _{L^{2}}  &  \leq Cc_{j}^{2}\left(  t\right)
\left(  H^{2}+U_{s}\left(  U_{s}+H\right)  \right)  . \label{R2}%
\end{align}
We get that%
\begin{align}
T_{1}  &  =\frac{\varepsilon}{2}\int\operatorname{div}V\left(  \eta_{j}%
^{2}+V_{j}^{2}\right)  +\varepsilon\int R_{1j}\eta_{j}+\varepsilon\int
R_{2j}V_{j}\nonumber\\
&  \leq\varepsilon\left\Vert \operatorname{div}V\right\Vert _{L^{\infty}%
}\left(  \left\Vert \eta_{j}\right\Vert _{L^{2}}^{2}+\left\Vert V_{j}%
\right\Vert _{L^{2}}^{2}\right)  +\varepsilon Cc_{j}^{3}\left(  t\right)
U_{j}\left(  H^{2}+U_{s}\left(  U_{s}+H\right)  \right) \nonumber\\
&  \leq C\varepsilon\left(  U_{j}^{2}U_{s}+c_{j}^{3}\left(  t\right)
U_{j}\left(  H^{2}+HU_{s}+U_{s}^{2}\right)  \right)  . \label{eq9'}%
\end{align}
We turn our attention towards $T_{2}$:%
\begin{align*}
T_{2}  &  =c\varepsilon^{2}\int(\nabla^{2}\eta_{j}W)\nabla\eta_{j}%
+c\varepsilon^{2}\int(\nabla^{2}\eta_{j}V)\nabla\eta_{j}+c\varepsilon^{2}%
\int\operatorname{div}V_{j}\nabla\eta\nabla\eta_{j}+c\varepsilon^{2}\int\nabla
V_{j}\nabla\eta\nabla\eta_{j}\\
&  +c\varepsilon^{2}\int\Delta_{j}\left(  \nabla W\nabla\eta\right)
\nabla\eta_{j}-c\varepsilon^{2}\int R_{3j}\nabla\eta_{j}.
\end{align*}
According to Proposition \ref{propA} \ and Proposition \ref{comut} we get
that\bigskip%
\begin{align}
&  \left\Vert \Delta_{j}\left(  \nabla W\nabla\eta\right)  \right\Vert
_{L^{2}}+\left\Vert R_{3j}\right\Vert _{L^{2}}\nonumber\\
&  \leq Cc_{j}^{4}\left(  t\right)  \left(  \left\Vert \nabla W\right\Vert
_{B_{2,r}^{s}}\left\Vert \nabla\eta\right\Vert _{B_{2,r}^{s}}+\left\Vert
\nabla W\right\Vert _{B_{2,r}^{s-1}}\left\Vert \nabla\eta\right\Vert
_{B_{2,r}^{s}}\right. \nonumber\\
&  \left.  \text{ \ \ \ \ \ \ \ \ \ \ \ \ \ \ \ \ }+\left\Vert \nabla
V\right\Vert _{B_{2,r}^{s-1}}\left\Vert \nabla\eta\right\Vert _{B_{2,r}^{s}%
}+\left\Vert \nabla\eta\right\Vert _{B_{2,r}^{s-1}}\left\Vert \nabla
V\right\Vert _{B_{2,r}^{s}}+\left\Vert \nabla^{2}\eta\right\Vert
_{B_{2,r}^{s-1}}\left\Vert V\right\Vert _{B_{2,r}^{s}}\right) \nonumber\\
&  \leq Cc_{j}^{4}\left(  t\right)  \left(  \left\Vert \nabla W\right\Vert
_{B_{2,r}^{s}}\left\Vert \nabla\eta\right\Vert _{B_{2,r}^{s}}+\left\Vert
W\right\Vert _{B_{2,r}^{s}}\left\Vert \nabla\eta\right\Vert _{B_{2,r}^{s}%
}\right. \nonumber\\
&  \left.  ~~~~~\ ~~~~~~~~~~~+\left\Vert \eta\right\Vert _{B_{2,r}^{s}%
}\left\Vert \nabla V\right\Vert _{B_{2,r}^{s}}+\left\Vert \nabla
\eta\right\Vert _{B_{2,r}^{s}}\left\Vert V\right\Vert _{B_{2,r}^{s}}\right)  .
\label{part}%
\end{align}
We observe that it is here that we need the restriction \eqref{abcd1} on the
parameters. Indeed if $a=d=0$ and $c<0$, then we only have $V\in B_{2,r}^{s}$
hence $U_{s}$ cannot control $\left\Vert \nabla V\right\Vert _{B_{2,r}^{s}}$.
Let us observe that:%
\[
\int(\nabla^{2}\eta_{j}W)\nabla\eta_{j}=-\frac{1}{2}\int\left\vert \nabla
\eta_{j}\right\vert ^{2}\operatorname{div}W=0
\]
and%
\[
\int(\nabla^{2}\eta_{j}V)\nabla\eta_{j}=-\frac{1}{2}\int\left\vert \nabla
\eta_{j}\right\vert ^{2}\operatorname{div}V.
\]
We infer that:%
\begin{align}
&  T_{2}=c\varepsilon^{2}\int\operatorname{div}V_{j}\nabla\eta\nabla\eta
_{j}+c\varepsilon^{2}\int(\nabla V_{j}\nabla\eta)\nabla\eta_{j}+c\varepsilon
^{2}\int\Delta_{j}\left(  \nabla W\nabla\eta\right)  \nabla\eta_{j}%
-c\varepsilon^{2}\int R_{3j}\nabla\eta_{j}\nonumber\\
&  \leq\frac{-c}{2}\varepsilon^{2}\left\{  \left\Vert \operatorname{div}%
V\right\Vert _{L^{\infty}}\left\Vert \nabla\eta_{j}\right\Vert _{L^{2}}%
^{2}+2\left\Vert \nabla\eta\right\Vert _{L^{\infty}}\left\Vert \nabla\eta
_{j}\right\Vert _{L^{2}}\left\Vert \nabla V_{j}\right\Vert _{L^{2}}\right.
\nonumber\\
&  ~~~~~~~~~~~~~~~~~+Cc_{j}^{4}\left(  t\right)  \left\Vert \nabla\eta
_{j}\right\Vert _{L^{2}}\left(  \left\Vert \nabla W\right\Vert _{B_{2,r}^{s}%
}\left\Vert \nabla\eta\right\Vert _{B_{2,r}^{s}}+\left\Vert W\right\Vert
_{B_{2,r}^{s}}\left\Vert \nabla\eta\right\Vert _{B_{2,r}^{s}}\right.
\nonumber\\
&  \left.  \left.
~~~~~~~~~~~~~~~~~~\ \ \ \ \ \ \ ~\ \ \ ~~\ \ \ ~~~\ \ \ ~~~\ \ +\left\Vert
\eta\right\Vert _{B_{2,r}^{s}}\left\Vert \nabla V\right\Vert _{B_{2,r}^{s}%
}+\left\Vert \nabla\eta\right\Vert _{B_{2,r}^{s}}\left\Vert V\right\Vert
_{B_{2,r}^{s}}\right)  \right\} \nonumber\\
&  \leq\varepsilon C\max\left(  \frac{-c}{b-c},\frac{-c}{d-a}\right)  \left(
U_{j}^{2}U_{s}+c_{j}^{4}\left(  t\right)  U_{j}HU_{s}\right)  . \label{eq10}%
\end{align}
\ We let%
\begin{equation}
C_{abcd}^{2}=\max\left(  \frac{-c}{b-c},\frac{-c}{d-a}\right)  . \label{2}%
\end{equation}
Let us estimate $T_{3}$. As above, owing to Proposition \ref{comut} we may
bound $R_{4j}$ in the following manner:
\begin{align}
-a\varepsilon^{2}\left\Vert R_{4j}\right\Vert _{L^{2}}  &  \leq-a\varepsilon
^{2}Cc_{j}^{5}\left(  t\right)  \left(  \left\Vert \nabla W\right\Vert
_{B_{2,r}^{s-1}}\left\Vert \nabla V\right\Vert _{B_{2,r}^{s}}+\left\Vert
\nabla V\right\Vert _{B_{2,r}^{s-1}}\left\Vert \nabla V\right\Vert
_{B_{2,r}^{s}}+\left\Vert \nabla\nabla^{t}V\right\Vert _{B_{2,r}^{s-1}%
}\left\Vert V\right\Vert _{B_{2,r}^{s}}\right) \nonumber\\
&  \leq-a\varepsilon^{2}Cc_{j}^{5}\left(  t\right)  \left(  \left\Vert
W\right\Vert _{B_{2,r}^{s}}\left\Vert \nabla V\right\Vert _{B_{2,r}^{s}%
}+\left\Vert V\right\Vert _{B_{2,r}^{s}}\left\Vert \nabla V\right\Vert
_{B_{2,r}^{s}}\right) \nonumber\\
&  \leq C\varepsilon^{\frac{3}{2}}\frac{-a}{\sqrt{d-a}}c_{j}^{5}\left(
t\right)  U_{s}\left(  H+U_{s}\right)  . \label{estim2}%
\end{align}
Also, we can write due to Proposition \ref{propA}:%
\begin{gather*}
-a\varepsilon^{2}\left(  \frac{1}{2}\left\Vert \Delta_{j}\nabla^{2}\left\vert
W\right\vert ^{2}\right\Vert _{L^{2}}+\left\Vert \Delta_{j}\left(  \nabla
W\nabla^{t}V\right)  \right\Vert _{L^{2}}+\left\Vert \Delta_{j}\left(
\nabla^{2}W:V\right)  \right\Vert _{L^{2}}+\left\Vert \Delta_{j}\left(  \nabla
V\nabla^{t}W\right)  \right\Vert _{L^{2}}\right) \\
\leq-a\varepsilon^{2}Cc_{j}^{6}\left(  t\right)  \left(  \left\Vert \nabla
^{2}W\right\Vert _{B_{2,r}^{s}}\left\Vert W\right\Vert _{B_{2,r}^{s}%
}+\left\Vert \nabla W\right\Vert _{B_{2,r}^{s}}^{2}+\left\Vert \nabla
W\right\Vert _{B_{2,r}^{s}}\left\Vert \nabla V\right\Vert _{B_{2,r}^{s}%
}\right. \\
\left.  ~~~~~~~~~+\left\Vert \nabla^{2}W\right\Vert _{B_{2,r}^{s}}\left\Vert
V\right\Vert _{B_{2,r}^{s}}+\left\Vert \nabla W\right\Vert _{B_{2,r}^{s}%
}\left\Vert \nabla V\right\Vert _{B_{2,r}^{s}}\right) \\
\leq C\varepsilon^{\frac{3}{2}}\max\left(  -a,\frac{-a}{\sqrt{d-a}}\right)
c_{j}^{6}\left(  t\right)  \left(  H^{2}+HU_{s}\right)  .
\end{gather*}
Then, using the integration by parts identity \eqref{identity} we get that:%
\begin{gather}
T_{3}=a\varepsilon^{2}\int\left(  \nabla^{2}V_{j}:W\right)  :\nabla
V_{j}+a\varepsilon^{2}\int\left(  \nabla^{2}V_{j}:V\right)  :\nabla
V_{j}\nonumber\\
+a\varepsilon^{2}\int\nabla V_{j}\nabla^{t}V:\nabla V_{j}-a\varepsilon^{2}\int
R_{4j}:\nabla V_{j}\nonumber\\
+a\varepsilon^{2}\int\left(  \frac{1}{2}\Delta_{j}\nabla^{2}\left\vert
W\right\vert ^{2}+\Delta_{j}\left(  \nabla W\nabla^{t}V\right)  +\Delta
_{j}\left(  \nabla^{2}W:V\right)  +\Delta_{j}\left(  \nabla V\nabla
^{t}W\right)  \right)  :\nabla V_{j}\nonumber\\
\leq-a\varepsilon^{2}\left(  \frac{1}{2}\left\Vert \operatorname{div}%
V\right\Vert _{L^{\infty}}\left\Vert \nabla V_{j}\right\Vert _{L^{2}}%
^{2}+\left\Vert \nabla^{t}V\right\Vert _{L^{\infty}}\left\Vert \nabla
V_{j}\right\Vert _{L^{2}}^{2}\right)  +\nonumber\\
+C\varepsilon\max\left(  \frac{-a}{\sqrt{d-a}},\frac{-a}{d-a}\right)  \left(
c_{j}^{5}\left(  t\right)  \varepsilon^{-\frac{1}{2}}U_{j}U_{s}\left(
H+U_{s}\right)  +c_{j}^{6}\left(  t\right)  \varepsilon^{-\frac{1}{2}}%
U_{j}\left(  H^{2}+HU_{s}\right)  \right) \nonumber\\
\leq\varepsilon C\max\left(  \frac{-a}{\sqrt{d-a}},\frac{-a}{d-a}\right)
\left(  U_{j}^{2}U_{s}+c_{j}^{7}\left(  t\right)  U_{j}\left(  H^{2}%
+HU_{s}+U_{s}^{2}\right)  \right)  \label{eq11}%
\end{gather}
and as before, let us denote by%
\begin{equation}
C_{abcd}^{3}=\max\left(  \frac{-a}{\sqrt{d-a}},\frac{-a}{d-a}\right)  .
\label{3}%
\end{equation}
Finally let us turn our attention towards $T_{4}$. Let us write:%
\[
T_{4}=G_{4}+B_{4},
\]
with%
\[
B_{4}=\frac{\varepsilon}{2}\int\partial_{t}\eta\left\vert V_{j}\right\vert
^{2}+\varepsilon^{2}\frac{d}{2}\int\partial_{t}\eta\nabla V_{j}:\nabla
V_{j}-\varepsilon^{2}d\left\langle \partial_{t}\nabla V_{j},\nabla\eta\otimes
V_{j}\right\rangle _{L^{2}}%
\]
and observe that by integration by parts and \eqref{R2}%

\begin{align}
G_{4}  &  =-\varepsilon^{2}\int\left(  \nabla V_{j}W\right)  \eta
V_{j}-\varepsilon^{2}\int\left(  \nabla V_{j}V\right)  \eta V_{j}%
+\varepsilon^{2}\int\eta R_{2j}V_{j}\nonumber\\
&  \leq\frac{\varepsilon^{2}}{2}\int\left(  \left(  \operatorname{div}\left(
\eta W\right)  +\operatorname{div}\left(  \eta V\right)  \right)  \left\vert
V_{j}\right\vert ^{2}\right)  +C\varepsilon^{2}c_{j}^{2}\left(  t\right)
\left\Vert V_{j}\right\Vert _{L^{2}}\left\Vert \eta\right\Vert _{L^{\infty}%
}\left(  H^{2}+HU_{s}+U_{s}^{2}\right) \nonumber\\
&  \leq C\varepsilon\left(  U_{j}^{2}\left(  HU_{s}+U_{s}^{2}\right)
+c_{j}^{8}\left(  t\right)  U_{j}U_{s}\left(  H^{2}+HU_{s}+U_{s}^{2}\right)
\right)  . \label{estim3}%
\end{align}
Let us now analyze the term $B_{4}$. We begin with
\[
\text{ }\frac{\varepsilon}{2}\int\partial_{t}\eta\left\vert V_{j}\right\vert
^{2}\leq\varepsilon\left\Vert \partial_{t}\eta\right\Vert _{L^{\infty}%
}\left\Vert V_{j}\right\Vert _{L^{2}}^{2}\leq\frac{\varepsilon}{2}U_{j}%
^{2}\left\Vert \partial_{t}\eta\right\Vert _{L^{\infty}}
\]
We have that:%
\begin{align}
\left\Vert \partial_{t}\eta\right\Vert _{L^{\infty}}  &  =\left\Vert \left(
I-b\varepsilon\Delta\right)  ^{-1}\left[  \left(  I+a\varepsilon\Delta\right)
\operatorname{div}V+\varepsilon\operatorname{div}\left(  \eta(W+V)\right)
\right]  \right\Vert _{L^{\infty}}\nonumber\\
&  \leq\left\Vert \left(  I-b\varepsilon\Delta\right)  ^{-1}\left(
I+a\varepsilon\Delta\right)  \operatorname{div}V\right\Vert _{L^{\infty}%
}+\varepsilon\left\Vert \left(  I-b\varepsilon\Delta\right)  ^{-1}%
\operatorname{div}\left(  \eta(W+V)\right)  \right\Vert _{L^{\infty}%
}\nonumber\\
&  \leq\left\Vert \left(  I-b\varepsilon\Delta\right)  ^{-1}\left(
I+a\varepsilon\Delta\right)  \operatorname{div}V\right\Vert _{B_{2,1}%
^{\frac{n}{2}}}+\varepsilon\left\Vert W\eta\right\Vert _{B_{2,r}^{s}%
}+\varepsilon\left\Vert \eta V\right\Vert _{B_{2,r}^{s}} \label{eq12}%
\end{align}
If $b>0$ or $a=b=0$, then because the operator $\left(  I-b\varepsilon
\Delta\right)  ^{-1}\left(  I+a\varepsilon\Delta\right)  $ maps $L^{2}$ to
$L^{2}$ with norm independent of $\varepsilon$, we get that:%
\[
\left\Vert \partial_{t}\eta\right\Vert _{L^{\infty}}\leq C\max\left(
1,\frac{-a}{b}\right)  \left(  U_{s}+HU_{s}+U_{s}^{2}\right)
\]
If $b=0$, $a<0$ then $d>0$ and we see that we have
\begin{equation}
-a\varepsilon\left\Vert \operatorname{div}\Delta V\right\Vert _{B_{2,1}%
^{\frac{d}{2}}}\leq-a\varepsilon\left\Vert \nabla^{2}V\right\Vert
_{B_{2,1}^{\frac{d}{2}+1}}\leq-a\varepsilon\left\Vert \nabla^{2}V\right\Vert
_{B_{2,1}^{s}}\leq\sqrt{\frac{-a}{d}}U_{s}. \label{eq13}%
\end{equation}
We set%
\begin{equation}
C_{abcd}^{4}=\left\{
\begin{array}
[c]{c}%
\max\left(  1,\frac{-a}{b}\right)  \text{ if }b>0~\ \text{or }a=b=0,\\
\max\left(  1,\sqrt{\frac{-a}{d}}\right)  \text{ if }b=0\text{ and }a<0.
\end{array}
\right.  \label{4}%
\end{equation}
Thus, we get that%
\begin{equation}
\text{ }\frac{\varepsilon}{2}\int\partial_{t}\eta\left\vert V_{j}\right\vert
^{2}\leq\varepsilon CC_{abcd}^{4}U_{j}^{2}\left(  U_{s}+HU_{s}+U_{s}%
^{2}\right)  . \label{estimB40}%
\end{equation}

In a similar fashion we can treat the second term of $B_{4}$ thus obtaining%
\begin{equation}
\varepsilon^{2}\frac{d}{2}\int\partial_{t}\eta\nabla V_{j}:\nabla V_{j}\leq
C\frac{d}{d-a}C_{abcd}^{4}U_{j}^{2}\left(  U_{s}+HU_{s}+U_{s}^{2}\right)  ,
\label{estimB41}%
\end{equation}
and we set%
\begin{equation}
C_{abcd}^{5}=\frac{d}{d-a}C_{abcd}^{4}. \label{5}%
\end{equation}
\ \ \ \ \ \ \ Finally, the last term of $B_{4}$ is estimated as follows:
\[
-\varepsilon^{2}d\left\langle \partial_{t}\nabla V_{j},\nabla\eta\otimes
V_{j}\right\rangle _{L^{2}}\leq\varepsilon^{2}d\left\Vert \nabla
\eta\right\Vert _{L^{\infty}}\left\Vert \partial_{t}\nabla V_{j}\right\Vert
_{L^{2}}\left\Vert V_{j}\right\Vert _{L^{2}},
\]
and using
\[
\partial_{t}\nabla V_{j}=-\left(  I-\varepsilon d\Delta\right)  ^{-1}%
\nabla\left(  \nabla\eta_{j}+c\varepsilon\nabla\Delta\eta_{j}+\varepsilon
\nabla V_{j}W+\varepsilon\nabla V_{j}V-\varepsilon R_{2j}\right)
\]
we observe that%
\begin{align*}
\varepsilon d\left\Vert \partial_{t}\nabla V_{j}\right\Vert _{L^{2}}  &
\leq\sqrt{d}\left(  \left\Vert \eta_{j}\right\Vert _{L^{2}}-c\varepsilon
\left\Vert \nabla\eta_{j}\right\Vert _{L^{2}}+\varepsilon\left\Vert
W\right\Vert _{L^{\infty}}\left\Vert \nabla V_{j}\right\Vert _{L^{2}}\right.
\\
&  ~~~~~\ \ \ ~~\ \ \left.  +\varepsilon\left\Vert V\right\Vert _{L^{\infty}%
}\left\Vert \nabla V_{j}\right\Vert _{L^{2}}+\varepsilon\left\Vert
R_{2j}\right\Vert _{L^{2}}\right) \\
&  \leq CC_{abcd}^{6}\left(  U_{j}\left(  1+H+U_{s}\right)  +c_{j}^{2}\left(
t\right)  \left(  H^{2}+HU_{s}+U_{s}^{2}\right)  \right)  ,
\end{align*}
where
\begin{equation}
C_{abcd}^{6}=\max\left(  \sqrt{d},\frac{-c\sqrt{d}}{sgn\left(  d\right)
\sqrt{b-c}},\sqrt{\frac{d}{d-a}}\right)  , \label{6}%
\end{equation}
thus, we conclude that%
\begin{align}
-\varepsilon^{2}d\left\langle \partial_{t}\nabla V_{j},\nabla\eta\otimes
V_{j}\right\rangle _{L^{2}}  &  \leq\varepsilon CC_{abcd}^{6}U_{j}\left(
U_{j}\left(  U_{s}+HU_{s}+U_{s}^{2}\right)  +\right. \nonumber\\
&  \left.  +c_{j}^{2}\left(  t\right)  U_{s}\left(  H^{2}+HU_{s}+U_{s}%
^{2}\right)  \right)  . \label{estimB42}%
\end{align}

Combining estimates\ \eqref{estimB40}, \eqref{estimB41}, \eqref{estimB42} we
obtain:%
\begin{equation}
B_{4}\leq\varepsilon C\tilde{C}_{abcd}U_{j}\left\{  U_{j}(U_{s}+HU_{s}%
+U_{s}^{2})+c_{j}^{9}\left(  t\right)  U_{s}\left(  H^{2}+HU_{s}+U_{s}%
^{2}\right)  \right\}  , \label{estimB4}%
\end{equation}
where
\begin{equation}
\tilde{C}_{abcd}=\max_{i=\overline{1,6}}C_{abcd}^{i} \label{max}%
\end{equation}
is the maximum of all constants depending on $a,b,c,d$ that appear in
relations \eqref{1}, \eqref{2}, \eqref{3}, \eqref{4}, \eqref{5}, \eqref{6}.

Thus, supposing that we are not in the case:%
\[
a=d=0\text{ and }c<0,\text{ }b>0,\text{ }
\]
we are able to successfully bound the $T_{i}$'s.

Finally, after adding up the estimations \eqref{estim1}, \eqref{eq9'},
\eqref{eq10}, \eqref{eq11}, \eqref{estim3}, \eqref{estimB4} we obtain that
there exists a positive constant $C>0$ depending only on $n$ and $s$ but not
on $\varepsilon$ such that
\begin{gather}
\frac{1}{2}\partial_{t}\left(  \int\eta_{j}^{2}+\varepsilon\left(  b-c\right)
\left\vert \nabla\eta_{j}\right\vert ^{2}+\varepsilon^{2}(-c)b\left(
\nabla^{2}\eta_{j}:\nabla^{2}\eta_{j}\right)  \right) \nonumber\\
+\frac{1}{2}\partial_{t}\left(  \int\left(  1+\varepsilon\eta\right)
V_{j}^{2}+\varepsilon\left(  d-a+d\varepsilon\eta\right)  \left(  \nabla
V_{j}:\nabla V_{j}\right)  +\varepsilon^{2}(-a)d\left(  \nabla^{2}V_{j}%
:\nabla^{2}V_{j}\right)  \right) \nonumber\\
\leq\varepsilon C\tilde{C}_{abcd}U_{j}\left(  U_{j}\left(  U_{s}+HU_{s}%
+U_{s}^{2}\right)  +c_{j}(t)(1+U_{s})\left(  H^{2}+HU_{s}+U_{s}^{2}\right)
\right)  , \label{eq14'}%
\end{gather}
where $\left(  c_{j}\left(  t\right)  \right)  _{j}$ is a sequence with
$\left(  2^{js}c_{j}\left(  t\right)  \right)  _{j}\in\ell^{r}(\mathbb{Z)}$,
having norm $1$ and $\tilde{C}_{abcd}$ is defined in \eqref{max}. Actually for
the sake of simplicity, from now on we will not carry on the distinction
between constants that depend on the $a,b,c,d$ parameters and the ones
depending on the dimension $n$ and regularity index $s$.

\subsection{Another useful estimation}

At this point, working with the non-cavitation hypothesis:
\begin{equation}
1+\varepsilon\eta_{0}\left(  x\right)  \geq\alpha>0, \label{cond}%
\end{equation}
using estimate \eqref{eq14} and a bootstrap argument would be sufficient in
order to obtain long time existence result similar to the one obtained in
\cite{Saut1} (with some restriction on the value of $\varepsilon$ depending
upon the initial data and $\alpha$). However, proceeding in a slightly
different manner we can avoid the use of \eqref{cond} (although, all
physically relevant data will verify it as $1+\varepsilon\eta$ represents the
total height of the water over the flat bottom). We also stress out that the
estimates established in this section are only available for the parameters
$a$, $b$, $c$, $d$ verifying \eqref{1111} with the exception of the two cases
\eqref{abcd}. Let us investigate the following quantity:%
\[
\partial_{t}\left(  \frac{\varepsilon}{2}\left\Vert \eta\right\Vert
_{L^{\infty}}U_{j}^{2}\right)  =I_{1}+I_{2},
\]
where%
\[
\left\{
\begin{array}
[c]{c}%
2I_{1}=\varepsilon\left\Vert \eta\right\Vert _{L^{\infty}}\partial_{t}%
U_{j}^{2},\\
2I_{2}=\varepsilon U_{j}^{2}\partial_{t}\left\Vert \eta\right\Vert
_{L^{\infty}}.
\end{array}
\right.
\]
Owing to \eqref{eq7} we see that:%
\begin{align}
I_{1}  &  =\frac{1}{2}\varepsilon\left\Vert \eta\right\Vert _{L^{\infty}%
}\left(  T_{1}+T_{2}+T_{3}\right)  +\frac{1}{2}\varepsilon\left\Vert
\eta\right\Vert _{L^{\infty}}\left(  -\varepsilon\int\eta\eta_{j}%
\operatorname{div}V_{j}+c\varepsilon^{2}\int\eta\nabla\operatorname{div}%
V_{j}\nabla\eta_{j}\right) \nonumber\\
&  \leq C\varepsilon^{2}U_{j}U_{s}\left(  U_{j}U_{s}+c_{j}\left(  t\right)
\left(  H^{2}+HU_{s}+U_{s}^{2}\right)  \right)  +C\varepsilon^{2}\left\Vert
\eta\right\Vert _{L^{\infty}}^{2}\left\Vert \eta_{j}\right\Vert _{L^{2}%
}\left\Vert \operatorname{div}V_{j}\right\Vert _{L^{2}}\nonumber\\
&  +C\varepsilon^{3}\left\Vert \eta\right\Vert _{L^{\infty}}^{2}\left\Vert
\nabla\operatorname{div}V_{j}\right\Vert _{L^{2}}\left\Vert \nabla\eta
_{j}\right\Vert _{L^{2}}\nonumber\\
&  \leq C\varepsilon^{2}U_{j}U_{s}\left(  U_{j}U_{s}+c_{j}\left(  t\right)
\left(  H^{2}+HU_{s}+U_{s}^{2}\right)  \right)  +C\varepsilon^{\frac{3}{2}%
}U_{j}^{2}U_{s}^{2}+C\varepsilon^{\frac{3}{2}}U_{j}^{2}U_{s}^{2}\nonumber\\
&  \leq C\varepsilon U_{j}\left(  U_{j}U_{s}^{2}+c_{j}\left(  t\right)
U_{s}\left(  H^{2}+HU_{s}+U_{s}^{2}\right)  \right)  . \label{I1}%
\end{align}
\

\begin{remark}
\label{obs}Let us notice that the term $c\varepsilon^{2}\int\eta
\nabla\operatorname{div}V_{j}\nabla\eta_{j}$ raises some important issues. In
order to successfully estimate it (and thus in order to have the validity of
\eqref{I1}), we need the restriction \eqref{abcd} on the parameters $a,b,c,d$.
The idea is that when $c\not =0$, we must have\footnote{One of the unknown
functions $\eta,V$ must have at least $B_{2,r}^{s+2}$-regularity level while
the other one needs $B_{2,r}^{s+1}$-regularity level.}
\[
sgn\left(  b\right)  +sgn\left(  d\right)  -sgn\left(  a\right)  \geq2.
\]
In view of the fact that $b+d>0$, it transpires that we must exclude the
cases:%
\begin{align*}
a  &  =d=0,\text{ }b>0,\text{ }c<0\text{ and}\\
a  &  =b=0,\text{ }d>0,\text{ }c<0\text{.}%
\end{align*}
Also, it is worth announcing that establishing \eqref{I1} isn't the only place
where the restriction on the parameters is needed. As it will be soon
revealed, in order to obtain local existence of solutions we will again have
to bound $-c\varepsilon^{2}\int\eta\nabla\operatorname{div}V_{j}\nabla\eta
_{j}$ and the above considerations will have to apply.
\end{remark}

In order to handle $I_{2}$ we use the fact that the function $t\rightarrow
\left\Vert \eta\left(  t\right)  \right\Vert _{L^{\infty}}$ is locally
Lipschitz we get that a.e. $\partial_{t}\left\Vert \eta\left(  t\right)
\right\Vert _{L^{\infty}}$ exists and besides, a.e. in time we have%
\begin{align*}
\left\vert \partial_{t}\left\Vert \eta\left(  t\right)  \right\Vert
_{L^{\infty}}\right\vert  &  =\left\vert \lim_{s\rightarrow t}\frac{\left\Vert
\eta\left(  t\right)  \right\Vert _{L^{\infty}}-\left\Vert \eta\left(
s\right)  \right\Vert _{L^{\infty}}}{t-s}\right\vert \leq\lim_{s\rightarrow
t}\left\Vert \frac{\eta\left(  t\right)  -\eta\left(  s\right)  }%
{t-s}\right\Vert _{L^{\infty}}\\
&  \leq\left\Vert \partial_{t}\eta\left(  t\right)  \right\Vert _{L^{\infty}%
}\leq\left\Vert \partial_{t}\eta\left(  t\right)  \right\Vert _{B_{2,1}%
^{\frac{n}{2}}}.
\end{align*}
Thus, owing to \eqref{eq12}, \eqref{eq13} we get that:%
\begin{equation}
I_{2}\leq\frac{1}{2}\varepsilon U_{j}^{2}\partial_{t}\left\Vert \eta
\right\Vert _{L^{\infty}}\leq C\varepsilon U_{j}^{2}\left(  U_{s}+HU_{s}%
+U_{s}^{2}\right)  . \label{I2}%
\end{equation}
Combining \eqref{I1} and \eqref{I2} we get that%
\begin{equation}
\partial_{t}\left(  \frac{\varepsilon}{2}\left\Vert \eta\right\Vert
_{L^{\infty}}U_{j}^{2}\right)  \leq C\varepsilon U_{j}\left(  U_{j}\left(
U_{s}+HU_{s}+U_{s}^{2}\right)  +c_{j}\left(  t\right)  U_{s}\left(
H^{2}+HU_{s}+U_{s}^{2}\right)  \right)  . \label{repar}%
\end{equation}
Finally let us observe that by adding \eqref{eq14} with \eqref{repar} we
obtain that:%
\begin{equation}
\partial_{t}N_{j}^{2}\leq C\varepsilon U_{j}\left(  U_{j}\left(  U_{s}%
+HU_{s}+U_{s}^{2}\right)  +c_{j}\left(  t\right)  \left(  1+U_{s}\right)
\left(  H^{2}+HU_{s}+U_{s}^{2}\right)  \right)  \label{eqnou}%
\end{equation}
where $N_{j}$ is the quantity defined in \eqref{Nj}.

\subsection{The case $b>0$, $c<0$ and $a=d=0$\label{Section}}

As pointed out in Section \ref{Sec} , we are not able to establish
\eqref{eq14} for the case \eqref{abcd1}. However, if we proceed slightly
different we can repair this inconvenience. Let us give some details of this
aspect. As we have seen, the problem comes when estimating $\left[
\eta,\Delta_{j}\right]  \nabla\operatorname{div}V$ (see \eqref{part}). In
order to bypass this problem, let us rewrite equation \eqref{eq5'} in the
following manner%
\begin{gather}
\frac{-c\varepsilon}{2}\partial_{t}\left(  \int\left(  \left\vert \nabla
\eta_{j}\right\vert ^{2}+\varepsilon b\nabla^{2}\eta_{j}:\nabla^{2}\eta
_{j}\right)  \right)  -c\varepsilon\int\nabla\operatorname{div}V_{j}\nabla
\eta_{j}\label{altern1}\\
-c\varepsilon^{2}\int\Delta_{j}\left(  \eta\nabla\operatorname{div}V\right)
\nabla\eta_{j}=T_{5}\nonumber
\end{gather}
where%
\begin{gather*}
T_{5}=c\varepsilon^{2}\int\left(  \nabla^{2}\eta_{j}W:\nabla\eta_{j}\right)
+c\varepsilon^{2}\int\left(  \nabla^{2}\eta_{j}V:\nabla\eta_{j}\right) \\
+c\varepsilon^{2}\left(  \int\operatorname{div}V_{j}\nabla\eta\nabla\eta
_{j}+\int\nabla V_{j}\nabla\eta\nabla\eta_{j}\right)  +c\varepsilon^{2}%
\int\Delta_{j}\left(  \nabla W\nabla\eta\right)  \nabla\eta_{j}-c\varepsilon
^{2}\int\tilde{R}_{j3}\nabla\eta_{j}%
\end{gather*}
and%
\[
\tilde{R}_{j3}=\left[  W,\Delta_{j}\right]  \nabla^{2}\eta+\left[
V,\Delta_{j}\right]  \nabla^{2}\eta+[\nabla\eta,\Delta_{j}]\operatorname{div}%
V+\left[  \nabla\eta,\Delta_{j}\right]  \nabla V.
\]
We add \eqref{altern1} with \eqref{eq4'} and \eqref{eq8} in order to obtain%
\begin{gather*}
\frac{1}{2}\partial_{t}\left(  \int\eta_{j}^{2}+\varepsilon(b-c)\left\vert
\nabla\eta_{j}\right\vert ^{2}+\varepsilon^{2}b\left(  -c\right)  \nabla
^{2}\eta_{j}:\nabla^{2}\eta_{j}+\left(  1+\varepsilon\eta\right)  \left\vert
V_{j}\right\vert ^{2}\right) \\
+c\varepsilon^{2}\int\eta\nabla\Delta\eta_{j}V_{j}-c\varepsilon^{2}\int
\Delta_{j}\left(  \eta\nabla\operatorname{div}V\right)  \nabla\eta
_{j}=\varepsilon\int\nabla\eta\eta_{j}V_{j}+T_{1}+T_{4}+T_{5}.
\end{gather*}
\ Let us write that:%
\begin{align*}
c\varepsilon^{2}\int\eta\nabla\Delta\eta_{j}V_{j}  &  =-c\varepsilon^{2}%
\int\nabla\eta\Delta\eta_{j}V_{j}-c\varepsilon^{2}\int\eta\Delta\eta
_{j}\operatorname{div}V_{j}\\
&  =-c\varepsilon^{2}\int\nabla\eta\Delta\eta_{j}V_{j}-c\varepsilon^{2}%
\int\Delta\eta_{j}[\eta,\Delta_{j}]\operatorname{div}V-c\varepsilon^{2}%
\int\Delta\eta_{j}\Delta_{j}\left(  \eta\operatorname{div}V\right) \\
&  =-c\varepsilon^{2}\int\nabla\eta\Delta\eta_{j}V_{j}-c\varepsilon^{2}%
\int\Delta\eta_{j}[\eta,\Delta_{j}]\operatorname{div}V+c\varepsilon^{2}%
\int\nabla\eta_{j}\Delta_{j}\left(  \nabla\eta\operatorname{div}V\right) \\
&  +c\varepsilon^{2}\int\nabla\eta_{j}\Delta_{j}\left(  \eta\nabla
\operatorname{div}V\right) \\
&  =-c\varepsilon^{2}\int\nabla\eta\Delta\eta_{j}V_{j}-c\varepsilon^{2}%
\int\Delta\eta_{j}[\eta,\Delta_{j}]\operatorname{div}V+c\varepsilon^{2}%
\int\nabla\eta_{j}\nabla\eta\operatorname{div}V_{j}\\
&  +c\varepsilon^{2}\int\nabla\eta_{j}[\Delta_{j},\nabla\eta
]\operatorname{div}V+c\varepsilon^{2}\int\nabla\eta\Delta_{j}\left(
\eta\nabla\operatorname{div}V\right)  .
\end{align*}
Thus, we get:
\begin{align}
&  c\varepsilon^{2}\int\eta\nabla\Delta\eta_{j}V_{j}-c\varepsilon^{2}%
\int\nabla\eta\Delta_{j}\left(  \eta\nabla\operatorname{div}V\right)
\nonumber\\
&  =-c\varepsilon^{2}\int\nabla\eta\Delta\eta_{j}V_{j}-c\varepsilon^{2}%
\int\Delta\eta_{j}[\eta,\Delta_{j}]\operatorname{div}V+c\varepsilon^{2}%
\int\nabla\eta_{j}\nabla\eta\operatorname{div}V_{j}\nonumber\\
&  +c\varepsilon^{2}\int\nabla\eta_{j}[\Delta_{j},\nabla\eta
]\operatorname{div}V\nonumber\\
&  =-c\varepsilon^{2}\int\nabla\eta\Delta\eta_{j}V_{j}-c\varepsilon^{2}%
\int\Delta\eta_{j}[\eta,\Delta_{j}]\operatorname{div}V-c\varepsilon^{2}%
\int\nabla^{2}\eta_{j}\nabla\eta V_{j}\nonumber\\
&  -c\varepsilon^{2}\int\nabla^{2}\eta\nabla\eta_{j}V_{j}+\int\nabla\eta
_{j}[\Delta_{j},\nabla\eta]\operatorname{div}V\nonumber\\
&  \leq-c\varepsilon^{2}\left(  \left\Vert \nabla\eta\right\Vert _{L^{\infty}%
}\left\Vert \Delta\eta_{j}\right\Vert _{L^{2}}\left\Vert V_{j}\right\Vert
_{L^{2}}+\left\Vert \Delta\eta_{j}\right\Vert _{L^{2}}c_{j}\left(  t\right)
\left\Vert \nabla\eta\right\Vert _{B_{2,r}^{s-1}}\left\Vert V\right\Vert
_{B_{2,r}^{s}}+\right. \nonumber\\
&  ~~~\ \ ~~\ \ \ \ \ \ \ ~+\left\Vert \nabla\eta\right\Vert _{L^{\infty}%
}\left\Vert \nabla^{2}\eta_{j}\right\Vert _{L^{2}}\left\Vert V_{j}\right\Vert
_{L^{2}}+\left\Vert \nabla^{2}\eta\right\Vert _{L^{\infty}}\left\Vert
\nabla\eta_{j}\right\Vert _{L^{2}}\left\Vert V_{j}\right\Vert _{L^{2}%
}+\nonumber\\
&  ~~\ \ \ \ \ \ \ \ ~\ ~\ ~\left.  +\left\Vert \nabla\eta_{j}\right\Vert
_{L^{2}}c_{j}\left(  t\right)  \left\Vert \nabla^{2}\eta\right\Vert
_{B_{2,r}^{s-1}}\left\Vert V\right\Vert _{B_{2,r}^{s}}\right) \nonumber\\
&  \leq-cC\varepsilon\left(  U_{j}^{2}U_{s}+c_{j}\left(  t\right)  U_{j}%
U_{s}^{2}\right)  . \label{estim5}%
\end{align}
A similar reasoning as in \eqref{part} together with estimates \eqref{eq9'},
\eqref{estim3}, \eqref{estimB4} and \eqref{estim5} gives us :%
\begin{gather}
\frac{1}{2}\partial_{t}\left(  \int\eta_{j}^{2}+\varepsilon\left(  b-c\right)
\left\vert \nabla\eta_{j}\right\vert ^{2}+\varepsilon^{2}(-c)b\left(
\nabla^{2}\eta_{j}:\nabla^{2}\eta_{j}\right)  +\left(  1+\varepsilon
\eta\right)  V_{j}^{2}\right) \nonumber\\
\leq C\varepsilon U_{j}\left(  U_{j}\left(  U_{s}+HU_{s}+U_{s}^{2}\right)
+c_{j}(t)(1+U_{s})\left(  H^{2}+U_{s}H+U_{s}^{2}\right)  \right)  .
\label{estims}%
\end{gather}

\section{Existence and uniqueness\label{S3}}

We begin by defining what kind of solutions we are looking for:%
\[
\left\{
\begin{array}
[c]{l}%
\left(  I-\varepsilon b\Delta\right)  \partial_{t}\eta+\operatorname{div}%
V+a\varepsilon\operatorname{div}\Delta V+\varepsilon W\nabla\eta
+\varepsilon\operatorname{div}\left(  \eta V\right)  =0,\\
\left(  I-\varepsilon d\Delta\right)  \partial_{t}V+\nabla\eta+c\varepsilon
\nabla\Delta\eta+\varepsilon\frac{1}{2}\nabla\left\vert W\right\vert
^{2}+\varepsilon\nabla WV+\varepsilon\nabla VW+\varepsilon\frac{1}{2}%
\nabla\left\vert V\right\vert ^{2}=0,\\
\eta_{|t=0}=\eta_{0}\text{ , }V_{|t=0}=V_{0}.
\end{array}
\right.
\]

\begin{definition}
\label{Definitie}Let us consider a positive time $T>0$ and $\left(  \eta
_{0},V_{0}\right)  \in L^{2}\times\left(  L^{2}\right)  ^{n}$ and $W\in\left(
H^{1}\left(  \mathbb{R}^{n}\right)  \right)  ^{n}$. A pair $\left(
\eta,V\right)  \in\mathcal{C}\left(  \left[  0,T\right]  ,L^{2}\times\left(
L^{2}\right)  ^{n}\right)  $ is called a solution to \eqref{eq2} on $\left[
0,T\right]  $ if for any $\left(  \phi,\psi\right)  \in\mathcal{C}^{1}\left(
\left[  0,T\right]  ,\mathcal{S\times S}^{n}\right)  $ and for all
$t\in\left[  0,T\right]  $, the following identities hold true:%
\begin{align*}
&
{\displaystyle\int\limits_{0}^{t}}
\left\langle \eta,\left(  I-\varepsilon b\Delta\right)  \partial_{t}%
\phi\right\rangle _{L^{2}}+%
{\displaystyle\int\limits_{0}^{t}}
\left\langle V,(I+a\varepsilon\Delta)\nabla\phi\right\rangle _{L^{2}%
}+\varepsilon%
{\displaystyle\int\limits_{0}^{t}}
\left\langle \eta,\nabla\left(  W\phi\right)  \right\rangle _{L^{2}%
}+\varepsilon%
{\displaystyle\int\limits_{0}^{t}}
\left\langle \eta V,\nabla\phi\right\rangle _{L^{2}}\\
&  =\left\langle \eta(t),\left(  I-\varepsilon b\Delta\right)  \phi
(t)\right\rangle _{L^{2}}-\left\langle \eta_{0},\left(  I-\varepsilon
b\Delta\right)  \phi(0)\right\rangle _{L^{2}}%
\end{align*}
and%
\begin{gather*}%
{\displaystyle\int\limits_{0}^{t}}
\left\langle V,\left(  I-\varepsilon d\Delta\right)  \partial_{t}%
\psi\right\rangle _{L^{2}}+%
{\displaystyle\int\limits_{0}^{t}}
\left\langle \eta,(I+c\varepsilon\Delta)\operatorname{div}\psi\right\rangle
_{L^{2}}+\varepsilon%
{\displaystyle\int\limits_{0}^{t}}
\left\langle \frac{\left\vert W\right\vert ^{2}}{2},\operatorname{div}%
\psi\right\rangle _{L^{2}}+\\
\varepsilon%
{\displaystyle\int\limits_{0}^{t}}
\left\langle V,\nabla^{t}W\psi\right\rangle _{L^{2}}++\varepsilon%
{\displaystyle\int\limits_{0}^{t}}
\left\langle V,\operatorname{div}\left(  \psi\otimes W\right)  \right\rangle
_{L^{2}}+\varepsilon%
{\displaystyle\int\limits_{0}^{t}}
\left\langle \frac{\left\vert V\right\vert ^{2}}{2},\operatorname{div}%
\psi\right\rangle _{L^{2}}\\
=\left\langle V(t),\left(  I-\varepsilon b\Delta\right)  \psi(t)\right\rangle
_{L^{2}}-\left\langle V_{0},\left(  I-\varepsilon b\Delta\right)
\psi(0)\right\rangle _{L^{2}}.
\end{gather*}

\end{definition}

Let us state now the following local existence and uniqueness theorem which
serves as an intermediary result for Theorem \ref{Teorema2}:

\begin{theorem}
\label{Teorema1}Let $a,b,c,d$ be chosen as in \eqref{1111} excluding the two
cases \eqref{abcd}, $b+d>0$, $r\in\left[  0,\infty\right]  $ and
$s\in\mathbb{R}$ such that:%
\begin{equation}
s>\frac{n}{2}+1\text{ or }s=\frac{n}{2}+1\text{ and }r=1. \label{rels}%
\end{equation}
Furthermore, let us consider $s_{1},$ $s_{2}$ and $s_{3}$ defined by relation
\eqref{relatie} and $W\in\left(  B_{2,r}^{s_{3}}\right)  ^{n}$.Then, for all
$\left(  \eta_{0},V_{0}\right)  \in B_{2,r}^{s_{1}}\times\left(
B_{2,r}^{s_{2}}\right)  ^{n}$ with $\operatorname{curl}V_{0}=0$, there exists
a positive $T>0$ and an unique solution
\begin{align*}
\left(  \eta,V\right)   &  \in\mathcal{C}\left(  [0,T],B_{2,r}^{s_{1}}%
\times\left(  B_{2,r}^{s_{2}}\right)  ^{n}\right)  \ \text{\ if }%
r<\infty\text{ or}\\
\left(  \eta,V\right)   &  \in L^{\infty}\left(  [0,T],B_{2,\infty}^{s_{1}%
}\times\left(  B_{2,\infty}^{s_{2}}\right)  ^{n}\right)  \cap%
{\displaystyle\bigcap\limits_{\beta>0}}
\mathcal{C}\left(  [0,T],B_{2,\infty}^{s_{1}-\beta}\times\left(  B_{2,\infty
}^{s_{2}-\beta}\right)  ^{n}\right)  \text{ if }r=\infty,
\end{align*}
of \ equation \eqref{eq2}. Moreover, if we denote by $T\left(  \eta_{0}%
,V_{0}\right)  $ the maximal time of existence then, if $T\left(  \eta
_{0},V_{0}\right)  <\infty$, we have that:%
\begin{align}
\lim_{t\rightarrow T\left(  \eta_{0},V_{0}\right)  }U_{s}\left(  t\right)   &
=\infty\text{ if }r<\infty,\label{exploz}\\
\limsup_{t\rightarrow T\left(  \eta_{0},V_{0}\right)  }U_{s}\left(  t\right)
&  =\infty\text{ if }r=\infty. \label{exploz2}%
\end{align}

\end{theorem}

Of course, the local existence result of the above theorem is not optimal. For
some particular choices of the $a$,$b$,$c$,$d$ parameters we can solve the
Cauchy problem for initial data in larger spaces, see for instance \cite{Anh},
\cite{Bona2}, \cite{Bona3}. However, the lower bound on the time of existence
is at most of order $O\left(  \varepsilon^{-\frac{1}{2}}\right)  $. Initially,
this is also the case of the solution constructed in Theorem \ref{Teorema1}.
As mentioned above this time scale is not satisfactory from a practical point
of view. However, as a by-product of the explosion criteria
\eqref{exploz}-\eqref{exploz2} of the solution of \eqref{eq2} and some refined
energy estimate, we can improve the lower bound of the $T\left(  \eta
_{0},V_{0}\right)  $ thus establishing $O\left(  \frac{1}{\varepsilon}\right)
$-long time existence.

Having established all the estimates that we need, let us proceed by proving
Theorem \ref{Teorema1}.

\begin{proof}
We will use the so called Friedrichs method. For all $m\in\mathbb{N}$, let us
consider $\mathbb{E}_{m}$ the low frequency cut-off operator defined by:%
\[
\mathbb{E}_{m}f=\mathcal{F}^{-1}\left(  \chi_{B\left(  0,m\right)  }\hat
{f}\right)  .
\]
We define the space
\[
L_{m}^{2}=\left\{  f\in L^{2}:\text{\textrm{Supp}}\hat{f}\subset B\left(
0,m\right)  \right\}
\]
which, endowed with the $\left\Vert \cdot\right\Vert _{L^{2}}$-norm is a
Banach space. Let us observe that due to Bernstein's lemma, all Sobolev norms
are equivalent on $L_{m}^{2}$. For all $m\in\mathbb{N}$, we consider the
following differential equation on $L_{m}^{2}$:%
\begin{equation}
\left\{
\begin{array}
[c]{l}%
\partial_{t}\eta=F_{m}\left(  \eta,V\right)  ,\\
\partial V=G_{m}\left(  \eta,V\right)  ,\\
\eta_{|t=0}=\mathbb{E}_{m}\eta_{0},\text{ }V_{|t=0}=\mathbb{E}_{m}V_{0},
\end{array}
\right.  \label{En2}%
\end{equation}
where $\left(  F_{m},G_{m}\right)  :L_{m}^{2}\times\left(  L_{m}^{2}\right)
^{n}\rightarrow L_{m}^{2}\times\left(  L_{m}^{2}\right)  ^{n}$ are defined by:%
\begin{align}
F_{m}\left(  \eta,V\right)   &  =-\mathbb{E}_{m}\left(  \left(  I-\varepsilon
b\Delta\right)  ^{-1}\left[  \left(  I+a\varepsilon\Delta\right)
\operatorname{div}V+\varepsilon W\nabla\eta+\varepsilon\operatorname{div}%
\left(  \eta V\right)  \right]  \right)  ,\label{Fm}\\
G_{m}\left(  \eta,V\right)   &  =-\mathbb{E}_{m}\left(  \left(  I-\varepsilon
d\Delta\right)  ^{-1}\left[  \left(  I+c\varepsilon\Delta\right)  \nabla
\eta+\frac{\varepsilon}{2}\nabla\left\vert V+W\right\vert ^{2}\right]
\right)  .\label{Gm}%
\end{align}
It transpires that due to the equivalence of the Sobolev norm, $\left(
F_{m},G_{m}\right)  $ is continuous and locally Lipschitz on $L_{m}^{2}%
\times\left(  L_{m}^{2}\right)  ^{n}$. Thus, the classical Picard theorem
ensures that there exists a nonnegative time $T_{m}>0$ and a unique solution
$\left(  \eta^{m},V^{m}\right)  :\mathcal{C}^{1}\left(  [0,T_{m}],L_{m}%
^{2}\times\left(  L_{m}^{2}\right)  ^{n}\right)  $. Moreover, if $T_{m}$ is
finite, then:%
\[
\lim_{t\rightarrow T_{m}}\left\Vert \left(  \eta^{m},V^{m}\right)  \right\Vert
_{L^{2}}=\infty.
\]
Another important aspect is that because of the property $\mathbb{E}_{m}%
^{2}=\mathbb{E}_{m}$ we get that the estimate obtained in \eqref{eq7} still
holds true for $\left(  \eta_{m},V_{m}\right)  $, namely
\begin{gather}
\frac{1}{2}\partial_{t}\left(  \int\left(  \eta_{j}^{m}\right)  ^{2}%
+\varepsilon\left(  b-c\right)  \left\vert \nabla\eta_{j}^{m}\right\vert
^{2}+\varepsilon^{2}(-c)b\left(  \nabla^{2}\eta_{j}^{m}:\nabla^{2}\eta_{j}%
^{m}\right)  \right)  \label{eq7m}\\
+\frac{1}{2}\partial_{t}\left(  \int\left\vert V_{j}^{m}\right\vert
^{2}+\varepsilon\left(  d-a\right)  \left(  \nabla V_{j}^{m}:\nabla V_{j}%
^{m}\right)  +\varepsilon^{2}(-a)d\left(  \nabla^{2}V_{j}^{m}:\nabla^{2}%
V_{j}^{m}\right)  \right)  \nonumber\\
+\varepsilon\int\eta^{m}\eta_{j}^{m}\operatorname{div}V_{j}^{m}-c\varepsilon
^{2}\int\eta^{m}\nabla\operatorname{div}V_{j}^{m}\nabla\eta_{j}^{m}%
=T_{1}+T_{2}+T_{3}\nonumber
\end{gather}
with $T_{1},T_{2},T_{3}$ defined as in relations \eqref{T1}, \eqref{T2},
\eqref{T3} but with $\left(  \eta^{m},V^{m}\right)  $ instead of $\left(
\eta,V\right)  $. Considering $U_{j}^{m}$ and $U_{s}^{m}$ the quantities
defined in \eqref{Uj} and \eqref{Us} with $\left(  \eta^{m},V^{m}\right)  $
instead of $\left(  \eta,V\right)  $. Also, for the $T_{i}$'s we dispose of
the estimates \eqref{eq9'}, \eqref{eq10} and \eqref{eq11} with $U_{j}^{m}$ and
$U_{s}^{m}$ instead of $U_{j}$ and $U_{s}$ and thus we gather that:%
\[
T_{1}+T_{2}+T_{3}\leq\varepsilon CU_{j}^{m}\left(  U_{j}^{m}U_{s}^{m}%
+c_{j}\left(  t\right)  \left(  H^{2}+HU_{s}^{m}+\left(  U_{s}^{m}\right)
^{2}\right)  \right)  .
\]
Also, let us notice that\footnote{Observe that due to \eqref{relatie} at least
one of $\eta,V$ has regularity level $B_{2,r}^{s+1}$ and thus, eventually by
an integration by parts we can obtain the announced estimate.}:%
\[
-\varepsilon\int\eta^{m}\eta_{j}^{m}\operatorname{div}V_{j}^{m}\leq
\varepsilon^{\frac{1}{2}}C\left(  U_{j}^{m}\right)  ^{2}U_{s}^{m}.
\]
Next, observe that\footnote{Again, because of \eqref{relatie} if $c\not =0$
then at least one of $\eta,V$ has regularity level $B_{2,r}^{s+2}$ while the
other one has regularity level $B_{2,r}^{s+1}$ see Remark \eqref{obs}. Thus,
eventually by an integration by parts, we get the announced result.}:%
\[
c\varepsilon^{2}\int\eta^{m}\nabla\operatorname{div}V_{j}^{m}\nabla\eta
_{j}^{m}\leq\varepsilon^{\frac{1}{2}}C\left(  U_{j}^{m}\right)  ^{2}U_{s}^{m}.
\]
We thus gather that:%
\[
\frac{1}{2}\partial_{t}\left(  U_{j}^{m}\right)  ^{2}\leq\varepsilon^{\frac
{1}{2}}CU_{j}^{m}\left(  U_{j}^{m}U_{s}^{m}+c_{j}\left(  t\right)  \left(
H^{2}+HU_{s}^{m}+\left(  U_{s}^{m}\right)  ^{2}\right)  \right)
\]
and by a Gronwall-type argument we obtain that for all $t\in\left[
0,T_{m}\right]  $:%
\[
U_{j}^{m}\left(  t\right)  \leq U_{j}^{m}\left(  0\right)  +\varepsilon
^{\frac{1}{2}}C%
{\displaystyle\int\limits_{0}^{t}}
\left(  U_{j}^{m}U_{s}^{m}+c_{j}\left(  \tau\right)  \left(  H^{2}+HU_{s}%
^{m}+\left(  U_{s}^{m}\right)  ^{2}\right)  \right)  d\tau,
\]
thus multiplying with $2^{js}$ and performing an $\ell^{r}(\mathbb{Z)}%
$-summation, owing to Minkowsky's theorem, we get that:%
\[
U_{s}^{m}\left(  t\right)  \leq U_{s}^{m}\left(  0\right)  +C\varepsilon
^{\frac{1}{2}}t+\varepsilon^{\frac{1}{2}}C\max\left(  1,H\right)
{\displaystyle\int\limits_{0}^{t}}
\left(  U_{s}^{m}+\left(  U_{s}^{m}\right)  ^{2}\right)  d\tau.
\]
We denote by $H^{\star}=$ $\max\left(  1,H\right)  $. Thus, Gronwall's lemma
along with the explosion criterion for ODE's gives us that:%
\[
\frac{\ln\left(  1+\frac{1}{U_{s}^{m}\left(  0\right)  }\right)  }%
{\varepsilon^{\frac{1}{2}}CH^{\star}}\leq T_{m}^{\star}%
\]
where $T_{m}^{\star}$ is the maximal time of existence for \eqref{En2}. Also,
due to the nature of the cut-off operator $\mathbb{E}_{m}$, it transpires
that:%
\[
U_{s}^{m}\left(  0\right)  \leq U_{s}\left(  0\right)  ,
\]
and consequently that:%
\begin{equation}
\frac{\ln\left(  1+\frac{1}{U_{s}\left(  0\right)  }\right)  }{\varepsilon
^{\frac{1}{2}}CH^{\star}}\leq T_{m}^{\star}.\label{timp}%
\end{equation}
In particular, for each time $T>0$ such that the above inequality is strictly
satisfied we obtain due to Gronwall's lemma that for all $t\in\left[
0,T\right]  $
\begin{equation}
U_{s}^{m}\left(  t\right)  \leq\frac{1}{e^{\ln\left(  1+\frac{1}{U_{s}\left(
0\right)  }\right)  -\varepsilon^{\frac{1}{2}}CH^{\star}T}-1}\label{unif}%
\end{equation}
and thus, the solution $\left(  \eta_{m},V_{m}\right)  $ is uniformly bounded
on $\left[  0,T\right]  $. Next, from relations \eqref{En2} and \eqref{Fm} we
get that $\left(  \partial_{t}\eta^{m}\right)  _{m\in\mathbb{N}}$ is uniformly
bounded on $\left[  0,T\right]  $ in $B_{2,r}^{s-1}$. Indeed, the first term
is%
\[
\left(  I-b\varepsilon\Delta\right)  ^{-1}\left(  I+a\varepsilon\Delta\right)
\operatorname{div}V^{m}\in B_{2,r}^{s_{4}}%
\]
with\footnote{remember that $b+d>0$}%
\[
s_{4}=s+2sgn(b)+sgn\left(  d\right)  +sgn(a)-1\geq s-1.
\]
In view of \eqref{rels}, $B_{2,r}^{s-1}$ is an algebra and thus the last two
terms of \eqref{Fm} are also at least in $B_{2,r}^{s-1}$. Thus, in view of the
uniform estimates \eqref{unif}, we get that $\left(  \partial_{t}\eta
^{m}\right)  _{m\in\mathbb{N}}$ is uniformly bounded on $\left[  0,T\right]  $
in $B_{2,r}^{s-1}$. Next, considering for all $p\in\mathbb{N}$ a smooth
function $\phi_{p}$ such that:%
\[
\left\{
\begin{array}
[c]{l}%
Supp\left(  \phi_{p}\right)  \subset B\left(  0,p+1\right)  \\
\phi_{p}=1\text{ on }B\left(  0,p\right)
\end{array}
\right.
\]
it follows, in view of Proposition \ref{compact} that for each $p\in
\mathbb{N}$, the sequence $\left(  \phi_{p}\eta^{m}\right)  _{m\in\mathbb{N}}$
is uniformly equicontiniuous on $\left[  0,T\right]  $ and that for all
$t\in\left[  0,T\right]  $, the set $\left\{  \phi_{p}\eta^{m}\left(
t\right)  :m\in\mathbb{N}\right\}  $ is relatively compact in $B_{2,r}^{s-1}$.
Thus, the Ascoli-Arzela Theorem combined with Cantor's diagonal process
provides us a subsequence of $\left(  \eta^{m}\right)  _{m\in\mathbb{N}}$
\footnote{Still denoted $\left(  \eta^{m}\right)  _{m\in\mathbb{N}}$ for the
sake of simplicity.
\par
{}}and a tempered distribution $\eta\in\mathcal{C}\left(  \left[  0,T\right]
,\mathcal{S}^{\prime}\right)  $ such that for all $\phi\in\mathcal{D}\left(
\mathbb{R}^{n}\right)  $:%
\[
\phi\eta^{m}\rightarrow\phi\eta\text{ in }\mathcal{C}\left(  \left[
0,T\right]  ,B_{2,r}^{s-1}\right)  \text{.}%
\]
Moreover, owing to Proposition \ref{PropBesov} and \eqref{unif} we get that
$\eta\in L_{T}^{\infty}\left(  B_{2,r}^{s_{1}}\right)  $ and using
interpolation we get that%
\[
\phi\eta^{m}\rightarrow\phi\eta\text{ in }\mathcal{C}\left(  \left[
0,T\right]  ,B_{2,r}^{s_{1}-\gamma}\right)
\]
for all $\gamma>0$. Of course, by the same argument one can get $V\in
L_{T}^{\infty}\left(  \left(  B_{2,r}^{s_{2}}\right)  ^{n}\right)  $ such that
for all $\psi\in\left(  \mathcal{D}\left(  \mathbb{R}^{n}\right)  \right)
^{n}$:%
\[
\psi V^{m}\rightarrow\phi V\text{ in }\mathcal{C}\left(  \left[  0,T\right]
,\left(  B_{2,r}^{s_{s}-\gamma}\right)  ^{n}\right)
\]
for all $\gamma>0.$ We claim that the properties enlisted above permit us to
pass to the limit when $m\rightarrow\infty$ in the equation verified $\eta
^{m}$ and $V^{m}$. By the Fatou property of Besov spaces we get that $\left(
\eta,V\right)  \in L_{t}^{\infty}(B_{2,r}^{s_{1}}\times\left(  B_{2,r}^{s_{2}%
}\right)  ^{n})$. It remains to verify that $\left(  \eta,V\right)  $ has the
announced regularity. Suppose that $r<\infty$. From $\eta$'s equation we see
that $\partial_{t}\eta\in L_{t}^{\infty}(B_{2,r}^{s-1})$ and thus, $\eta
\in\mathcal{C}\left(  [0,T],B_{2,r}^{s-1}\right)  $ which also implies that
$S_{j}\eta\in\mathcal{C}\left(  [0,T],B_{2,r}^{s_{1}}\right)  $ for all
$j\in\mathbb{Z}$. The conclusion follows as the sequence of $B_{2,r}^{s_{1}}%
-$valued functions $\left(  S_{j}\eta\right)  _{j\in\mathbb{Z}}$ tends
uniformly to $\eta$. Indeed, for all $\ell,j\in\mathbb{Z}$, we have
\[
\Delta_{\ell}(\eta-S_{j}\eta)=0\text{ if }\ell\leq j-2,
\]
therefore we have%
\begin{align*}
\left\Vert \eta-S_{j}\eta\right\Vert _{L_{t}^{\infty}\left(  B_{2,r}^{s_{1}%
}\right)  }  & \leq\left(  \sum_{\ell\geq j-1}2^{\ell rs_{1}}\left\Vert
\eta_{0}\right\Vert _{L^{2}}^{r}\right)  ^{\frac{1}{r}}\\
& +\varepsilon^{\frac{1}{2}}C\left(  \left\Vert U_{s}\right\Vert
_{L_{t}^{\infty}\left(  B_{2,r}^{s_{1}}\right)  }^{2}+H^{2}\right)
{\displaystyle\int\limits_{0}^{t}}
\left(  \sum_{\ell\geq j-1}2^{\ell rs_{1}}\left(  U_{j}(\tau)+c_{j}%
(\tau)\right)  ^{r}\right)  ^{\frac{1}{r}}%
\end{align*}
and as a consequence of Proposition \ref{PropBesov} and the dominated
convergence theorem we find that
\[
\lim_{j\rightarrow\infty}\left\Vert \eta-S_{j}\eta\right\Vert _{L_{t}^{\infty
}\left(  B_{2,r}^{s_{1}}\right)  }=0,
\]
which implies that $\eta\in\mathcal{C}\left(  [0,T],B_{2,r}^{s_{1}}\right)  $.
A similar argument shows that $V\in\mathcal{C}\left(  [0,T],\left(
B_{2,r}^{s_{2}}\right)  ^{n}\right)  $. When $r=\infty$, we know that for
every positive $\beta$, $B_{2,1}^{s_{1}-\beta}$ is continuously embedded in
$B_{2,\infty}^{s_{1}}$ and repeating the above argument permits us to conclude
that $\left(  \eta,V\right)  $ has the desired regularity. This completes the
proof of existence.
Uniqueness, is a consequence of the following stability estimate: let us
consider two solutions of \eqref{eq2}, $\left(  \eta^{1},V^{1}\right)  $,
$\left(  \eta^{2},V^{2}\right)  $ and observe that the difference
\[
\left(  \delta\eta,\delta V\right)  =\left(  \eta^{1}-\eta^{2},V^{1}%
-V^{2}\right)
\]
satisfies the following system:%
\begin{equation}
\left\{
\begin{array}
[c]{l}%
\left(  I-\varepsilon b\Delta\right)  \partial_{t}\delta\eta
+\operatorname{div}\delta V+a\varepsilon\operatorname{div}\Delta\delta
V+\varepsilon W\nabla\delta\eta+\varepsilon\operatorname{div}\left(
\delta\eta V^{1}\right)  +\varepsilon\operatorname{div}\left(  \eta^{2}\delta
V\right)  =0\\
\left(  I-\varepsilon d\Delta\right)  \partial_{t}\delta V+\nabla\delta
\eta+c\varepsilon\nabla\Delta\delta\eta+\varepsilon\nabla W(\delta
V)+\varepsilon\nabla(\delta V)W+\varepsilon\nabla(\delta V)V^{1}%
+\varepsilon\nabla V^{2}\delta V=0\\
\delta\eta_{|t=0}=0\text{, }\delta V_{|t=0}=0.
\end{array}
\right.  \label{Stab1}%
\end{equation}
We consider $U_{s}^{1}=U_{s}\left(  \eta^{1},V^{1}\right)  $ and $U_{s}%
^{2}=U_{s}\left(  \eta^{2},V^{2}\right)  $, see \eqref{Us}. For the sake of
simplicity we will prove stability estimates in the classical Sobolev space
$X=H^{r_{1}}\times\left(  H^{r_{2}}\right)  ^{n}$ with:%
\[
\left\{
\begin{array}
[c]{c}%
r_{1}=sgn\left(  b\right)  -sgn\left(  c\right)  ,\\
r_{2}=sgn\left(  d\right)  -sgn\left(  a\right)  .
\end{array}
\right.
\]
We endow $X$ with the norm:%
\begin{align*}
\left\Vert \eta,V\right\Vert _{X}^{2} &  =\left\Vert \eta\right\Vert _{L^{2}%
}^{2}+\varepsilon\left(  b-c\right)  \left\Vert \nabla\eta\right\Vert _{L^{2}%
}^{2}-\varepsilon^{2}bc\left\Vert \nabla^{2}\eta\right\Vert _{L^{2}}^{2}\\
&  +\left\Vert V\right\Vert _{L^{2}}^{2}+\varepsilon\left(  d-a\right)
\left\Vert \nabla V\right\Vert _{L^{2}}^{2}-\varepsilon^{2}da\left\Vert
\nabla^{2}V\right\Vert _{L^{2}}^{2}%
\end{align*}
Observe that due to the fact that $s_{1},s_{2}$ are chosen so as to satisfy
\eqref{relatie} with $s$ chosen as in \eqref{rels}, we have that $H^{r_{1}%
}\times\left(  H^{r_{2}}\right)  ^{n}$ is continuously embedded in
$B_{2,r}^{s_{1}}\times\left(  B_{2,r}^{s_{2}}\right)  ^{n}$. Multiplying the
first equation in \eqref{Stab1} with $\eta$, the second equation with $V$ and
adding up the results we get that:%
\begin{gather*}
\frac{1}{2}\partial_{t}\left(  \left\Vert \delta\eta\right\Vert _{L^{2}}%
^{2}+\varepsilon\left(  b-c\right)  \left\Vert \nabla\delta\eta\right\Vert
_{L^{2}}^{2}+\left\Vert \delta V\right\Vert _{L^{2}}^{2}+\varepsilon\left(
d-a\right)  \left\Vert \nabla\delta V\right\Vert _{L^{2}}^{2}\right)
+\varepsilon a\int\Delta\operatorname{div}\delta V\delta\eta\\
+\varepsilon c\int\Delta\nabla\delta\eta\delta V\leq C\sqrt{\varepsilon
}\left(  H+U_{s}^{1}+U_{s}^{2}\right)  \left\Vert \delta\eta,\delta
V\right\Vert _{X}^{2}.
\end{gather*}
Now, let us multiply the first equation in \eqref{Stab1} with $c\varepsilon
\Delta\eta$, the second equation with $a\varepsilon\Delta V$ and adding up the
results we obtain:%
\begin{align*}
&  \frac{1}{2}\partial_{t}\left(  -c\varepsilon\left\Vert \nabla\delta
\eta\right\Vert _{L^{2}}^{2}-\varepsilon^{2}bc\left\Vert \nabla^{2}\delta
\eta\right\Vert _{L^{2}}^{2}-\varepsilon a\left\Vert \nabla\delta V\right\Vert
_{L^{2}}^{2}-\varepsilon^{2}da\left\Vert \nabla^{2}\delta V\right\Vert
_{L^{2}}^{2}\right)  \\
&  +\varepsilon a\int\nabla\delta\eta\Delta\delta V+\varepsilon c\int
\operatorname{div}\delta V\Delta\delta\eta\\
&  \leq C\sqrt{\varepsilon}\left(  H+U_{s}^{1}+U_{s}^{2}\right)  \left\Vert
\delta\eta,\delta V\right\Vert _{X}^{2}.
\end{align*}
Thus, adding up the above relations gives us:%
\[
\partial_{t}\left\Vert \delta\eta,\delta V\right\Vert _{X}^{2}\leq
C\sqrt{\varepsilon}\left(  H+U_{s}^{1}+U_{s}^{2}\right)  \left\Vert \delta
\eta,\delta V\right\Vert _{X}^{2}.
\]
Hence, Gronwall's lemma ensures the desired result. \newline Proving the
blow-up criteria is classic. Let us suppose that $T\left(  \eta_{0}%
,V_{0}\right)  <\infty$ and that%
\[
\limsup_{t\rightarrow T\left(  \eta_{0},V_{0}\right)  }U_{s}\left(  t\right)
<\infty.
\]
Then, $U_{s}\left(  t\right)  $ remains bounded on $[0,T\left(  \eta_{0}%
,V_{0}\right)  )$ say:%
\[
U_{s}\left(  t\right)  \leq M,
\]
for all $t\in\lbrack0,T\left(  \eta_{0},V_{0}\right)  )$. We see that for any
$t_{0}\in$ $[0,T\left(  \eta_{0},V_{0}\right)  )$ we can construct, using the
same method as before a solution to \eqref{eq2} with initial data $\left(
\eta\left(  t_{0}\right)  ,V\left(  t_{0}\right)  \right)  $ on a time
interval that according to \eqref{timp} satisfies the following lower bound:%
\[
T^{new}-t_{0}\geq\frac{\ln\left(  1+\frac{1}{U_{s}\left(  t_{0}\right)
}\right)  }{\varepsilon^{\frac{1}{2}}CH^{\star}}\geq\frac{\ln\left(
1+\frac{1}{M}\right)  }{\varepsilon^{\frac{1}{2}}CH^{\star}}.
\]
Of course, choosing $t_{0}$ close enough to $T\left(  \eta_{0},V_{0}\right)  $
we can obtain $T^{new}>T\left(  \eta_{0},V_{0}\right)  $ such that gluing
together the new solution with the $\left(  \eta,V\right)  _{|\left[
0,t_{0}\right]  }$ and in view of the uniqueness we get a contradiction on the
maximality of $T\left(  \eta_{0},V_{0}\right)  $. This concludes the proof of
the announced result.
\end{proof}

\section{The proof of Theorem \ref{Teorema2}\label{iar}}

We are now in the position of establishing the announced long time existence
result. Actually, we prove long time existence and uniqueness in the framework
of the more general Besov space. More precisely, Theorem \ref{Teorema2} is
just a particular case of the following:

\begin{theorem}
\label{Teo}Let $a,b,c,d$ as in \eqref{1111} excluding the two
cases\ \eqref{abcd}, $b+d>0$. Let us take $r\in\left[  1,\infty\right]  $ and
$s$ such that%
\[
s>\frac{n}{2}+1\text{ or }s=\frac{n}{2}+1\text{ and }r=1,
\]
with $n\geq1$. Let us also consider $s_{1}$, $s_{2}$ and $s_{3}$ defined by
\eqref{relatie} and $W\in\left(  B_{2,r}^{s_{3}}\right)  ^{n}$. Then, we can
establish long time existence and uniqueness of solutions (see Definition
\ref{def}) for the equation \eqref{eq2} in $B_{2,r}^{s_{1}}\times\left(
B_{2,r}^{s_{2}}\right)  ^{n}$. Moreover, if we denote by $T\left(  \eta
_{0},V_{0}\right)  $, the maximal time of existence then there exists some
$T\in\lbrack0,T\left(  \eta_{0},V_{0}\right)  )$ which is bounded from below
by an $O\left(  \frac{1}{\varepsilon}\right)  $-order quantity and a function
$G:\mathbb{R\rightarrow R}$ such that for all $t\in\left[  0,T\right]  $ we
have:
\[
U_{s}\left(  \eta,V\right)  \leq G\left(  U_{s}\left(  \eta_{0},V_{0}\right)
\right)
\]
where $U_{s}(\eta,V)$ is defined in relation \eqref{Us}.
\end{theorem}

Let us consider the unique maximal solution $\left(  \eta,V\right)  $ of
\eqref{eq2} that we constructed in the proof of Theorem \ref{Teorema1}. Of
course, the estimate \eqref{eqnou} holds true for this solution thus we have
that:%
\begin{equation}
\partial_{t}N_{j}^{2}\leq C\varepsilon U_{j}\left(  U_{j}\left(  U_{s}%
+HU_{s}+U_{s}^{2}\right)  +c_{j}(t)(1+U_{s})\left(  H^{2}+U_{s}H+U_{s}%
^{2}\right)  \right)  . \label{eq15'}%
\end{equation}
We recall that $N_{j}$ is the following quantity:%

\begin{gather*}
N_{j}^{2}\left(  t\right)  =\int\left(  1+\varepsilon\left\Vert \eta
\right\Vert _{L^{\infty}}\right)  \eta_{j}^{2}+\varepsilon\left(  b-c\right)
\left(  1+\varepsilon\left\Vert \eta\right\Vert _{L^{\infty}}\right)
\left\vert \nabla\eta_{j}\right\vert ^{2}\\
+\int\varepsilon^{2}(-c)b\left(  1+\varepsilon\left\Vert \eta\right\Vert
_{L^{\infty}}\right)  \left(  \nabla^{2}\eta_{j}:\nabla^{2}\eta_{j}\right) \\
+\int\left(  1+\varepsilon\eta+\varepsilon\left\Vert \eta\right\Vert
_{L^{\infty}}\right)  V_{j}^{2}+\varepsilon\left(  d-a+d\varepsilon
\eta+d\varepsilon\left\Vert \eta\right\Vert _{L^{\infty}}\right)  \left(
\nabla V_{j}:\nabla V_{j}\right) \\
+\int\varepsilon^{2}(-a)d\left(  1+\varepsilon\left\Vert \eta\right\Vert
_{L^{\infty}}\right)  \left(  \nabla^{2}V_{j}:\nabla^{2}V_{j}\right)  .
\end{gather*}
and that we have:
\begin{equation}
U_{j}\left(  t\right)  \leq N_{j}\left(  t\right)  \leq\left(  1+2\varepsilon
\left\Vert \eta\left(  t\right)  \right\Vert _{L^{\infty}}\right)  ^{\frac
{1}{2}}U_{j}\left(  t\right)  \label{eq15}%
\end{equation}
thus we immediately get that%
\[
\partial_{t}N_{j}^{2}\leq C\varepsilon N_{j}\left(  U_{j}\left(  U_{s}%
+HU_{s}+U_{s}^{2}\right)  +c_{j}(t)(1+U_{s})\left(  H^{2}+U_{s}H+U_{s}%
^{2}\right)  \right)
\]
and by time integration we get that:%
\[
U_{j}\left(  t\right)  \leq N_{j}\left(  t\right)  \leq N_{j}\left(  0\right)
+C\varepsilon\int_{0}^{t}\left(  U_{j}\left(  U_{s}+HU_{s}+U_{s}^{2}\right)
+c_{j}(t)(1+U_{s})\left(  H^{2}+U_{s}H+U_{s}^{2}\right)  \right)  d\tau.
\]
Multiplying the last inequality with $2^{js}$ and performing a $\ell
^{r}(\mathbb{Z)}$ summation yields:%
\begin{align}
U_{s}\left(  t\right)   &  \leq N_{0}+C\varepsilon\int_{0}^{t}\left(
H^{2}+U_{s}\left(  H+H^{2}\right)  +U_{s}^{2}\left(  1+H\right)  +U_{s}%
^{3}\right)  d\tau,\nonumber\\
&  \leq N_{0}+C\varepsilon tH^{2}+\varepsilon C\left(  1+H+H^{2}\right)
\int_{0}^{t}\left(  U_{s}+U_{s}^{3}\right)  d\tau,\nonumber\\
&  \leq N_{0}+C\varepsilon tH^{2}+\varepsilon C\left(  1+H^{2}\right)
\int_{0}^{t}\left(  U_{s}+U_{s}^{3}\right)  d\tau, \label{eq16}%
\end{align}
for all $t\in\lbrack0,T\left(  \eta_{0},V_{0}\right)  )$ where
\[
N_{0}=\left\Vert \left(  2^{js}N_{j}\left(  0\right)  \right)  _{j\in
\mathbb{Z}}\right\Vert _{\ell^{r}(\mathbb{Z)}}\leq\left(  1+2\varepsilon
\left\Vert \eta_{0}\right\Vert _{L^{\infty}}\right)  ^{\frac{1}{2}}\left\Vert
\left(  \eta_{0},V_{0}\right)  \right\Vert _{B_{2,r}^{s}}.
\]
Having established \eqref{eq16} we are in the position of finding a lower
bound of order $O\left(  \frac{1}{\varepsilon}\right)  $ for the time of
existence. Indeed by Gronwall's lemma, and taking in account the explosion
criterion of Theorem \ref{Teorema1} we get that, supposing $T\left(  \eta
_{0},V_{0}\right)  $ is finite:
\begin{equation}
\limsup_{t\rightarrow T\left(  \eta_{0},V_{0}\right)  }U_{s}\left(  t\right)
=+\infty. \label{explosion}%
\end{equation}
Hence, via Gronwall's lemma, we have:%
\begin{equation}
\frac{1}{2}\ln\left(  1+\frac{1}{\left(  N_{0}+\varepsilon CH^{2}T\left(
\eta_{0},V_{0}\right)  \right)  ^{2}}\right)  \leq\varepsilon C\left(
1+H^{2}\right)  T\left(  \eta_{0},V_{0}\right)  . \label{concluzie}%
\end{equation}
If%
\[
\varepsilon CH^{2}T\left(  \eta_{0},V_{0}\right)  \geq N_{0}\text{ i.e.
}T\left(  \eta_{0},V_{0}\right)  \geq\frac{N_{0}}{\varepsilon CH^{2}},
\]
then, we have nothing else to prove. If this is not the case, the LHS member
of \eqref{concluzie} is larger then $\frac{1}{2}\ln\left(  1+\frac{1}%
{4N_{0}^{2}}\right)  $ and we get that:%
\[
T\left(  \eta_{0},V_{0}\right)  \geq\frac{\frac{1}{2}\ln\left(  1+\frac
{1}{4N_{0}^{2}}\right)  }{\varepsilon C\left(  1+H^{2}\right)  }.
\]
Thus, $T\left(  \eta_{0},V_{0}\right)  $ is always bounded from bellow by a
quantity of order $O\left(  \frac{1}{\varepsilon}\right)  $.

We now prove that on a $O\left(  \frac{1}{\varepsilon}\right)  $-order time
interval we dispose of uniform bounds for the solution of \eqref{eq2}. Let us
consider%
\[
T^{\star}=\sup\left\{  T\in\lbrack0,T\left(  \eta_{0},V_{0}\right)  ):\forall
t\in\left[  0,T\right]  ,\text{ }U_{s}\left(  t\right)  \leq2\left(
N_{0}+\varepsilon CtH^{2}\right)  \right\}  .
\]
Then, from \eqref{eq16}, we deduce that for all $t\leq T^{\star}$ we have%
\[
U_{s}\left(  t\right)  \leq N_{0}+\varepsilon CT^{\star}H^{2}+\varepsilon
C\left(  1+H^{2}\right)  \left(  1+4\left(  N_{0}+\varepsilon CT^{\star}%
H^{2}\right)  ^{2}\right)  \int_{0}^{t}U_{s}\left(  \tau\right)  d\tau
\]
and according to Gronwall's lemma we get that:%
\[
U_{s}\left(  t\right)  \leq\left(  N_{0}+\varepsilon CT^{\star}H^{2}\right)
\exp\left(  \varepsilon T^{\star}C\left(  1+H^{2}\right)  \left(  1+4\left(
N_{0}+\varepsilon CT^{\star}H^{2}\right)  ^{2}\right)  \right)  .
\]
Now, if there exists a $\beta\in\left(  0,2\right)  $ such that%
\[
\exp\left(  \varepsilon T^{\star}C\left(  1+H^{2}\right)  \left(  1+4\left(
N_{0}+\varepsilon CT^{\star}H^{2}\right)  ^{2}\right)  \right)  \leq\beta,
\]
then a continuity argument will lead us to the conclusion that $T^{\star
}=T\left(  \eta_{0},V_{0}\right)  $ which will imply that $T\left(  \eta
_{0},V_{0}\right)  =\infty$. Thus, for all $t\in\left[  0,\frac{N_{0}%
}{\varepsilon CH^{2}}\right]  $ we get that
\begin{equation}
U_{s}\left(  t\right)  \leq4N_{0}\leq\left(  1+2\varepsilon U_{s}\left(
\eta_{0},V_{0}\right)  \right)  ^{\frac{1}{2}}U_{s}\left(  \eta_{0}%
,V_{0}\right)  . \label{411}%
\end{equation}
Assume that%
\begin{equation}
\varepsilon T^{\star}C\left(  1+H^{2}\right)  \left(  1+4\left(
N_{0}+\varepsilon CT^{\star}H^{2}\right)  ^{2}\right)  \geq\ln2. \label{413}%
\end{equation}
Then, if
\begin{equation}
T^{\star}>\frac{N_{0}}{\varepsilon CH^{2}}, \label{412}%
\end{equation}
we get that \eqref{411} holds true for all $t\in\left[  0,\frac{N_{0}%
}{\varepsilon CH^{2}}\right]  $. If \eqref{412} doesn't hold i.e.%
\begin{equation}
T^{\star}\leq\frac{N_{0}}{\varepsilon CH^{2}} \label{414}%
\end{equation}
then, combining \eqref{413} with \eqref{414}, we get that:%
\[
T^{\star}\geq\frac{\ln2}{\varepsilon C\left(  1+H^{2}\right)  \left(
1+16N_{0}^{2}\right)  }%
\]
and thus, using the definition of $T^{\star}$, we get that for all
$t\in\left[  0,\frac{\ln2}{\varepsilon C\left(  1+H^{2}\right)  \left(
1+16N_{0}^{2}\right)  }\right]  $
\begin{align}
U_{s}\left(  t\right)   &  \leq2\left(  N_{0}+\frac{H^{2}\ln2}{\left(
1+H^{2}\right)  \left(  1+16N_{0}^{2}\right)  }\right)  ,\nonumber\\
&  \leq2\left(  N_{0}+\frac{\ln2}{\left(  1+16N_{0}^{2}\right)  }\right)
\label{415}\\
&  \leq2\left(  \left(  1+2\varepsilon U_{s}\left(  \eta_{0},V_{0}\right)
\right)  ^{\frac{1}{2}}U_{s}\left(  \eta_{0},V_{0}\right)  +\frac{\ln
2}{\left(  1+16U_{s}^{2}\left(  \eta_{0},V_{0}\right)  \right)  }\right)  .
\label{416}%
\end{align}
Thus, one can choose
\[
\left\{
\begin{array}
[c]{c}%
F(x)=\min\left\{  \frac{\left(  1+2\left\vert x\right\vert \right)  ^{\frac
{1}{2}}x}{CH^{2}},\frac{\ln2}{C\left(  1+H^{2}\right)  \left(  1+16x^{2}%
\right)  }\right\}  ,\\
G(x)=\max\left\{  2\left(  \left(  1+2\left\vert x\right\vert \right)
^{\frac{1}{2}}x+\frac{\ln2}{1+16x^{2}}\right)  ,4\left(  1+2\left\vert
x\right\vert \right)  ^{\frac{1}{2}}x\right\}  ,
\end{array}
\right.
\]
in order to conclude that $T^{\star}\geq F\left(  U_{0}\right)  $ and for all
$t\in\lbrack0,T^{\star}]$ we have:
\[
U_{s}(t)\leq G\left(  U_{0}\right)  .
\]

\begin{remark}
Theorem \ref{Teorema2} is just the restatement of the above result when $r=2$.
\end{remark}

\begin{remark}
Observe that the above arguments allow us to derive similar uniform bounds for
the quantities $\left\Vert \left(  2^{js}N_{j}\right)  _{j\in\mathbb{Z}%
}\right\Vert _{\ell^{r}(\mathbb{Z})}$.
\end{remark}

\section{Final remarks}

The results of this paper generalize a part of the long time existence results
for the $abcd$ Boussinesq systems that can be found in \cite{Saut1}. Theorem
\ref{Teo} has also a homogeneous counterpart i.e. one can replace the
nonhomogeneous Besov space $B_{2,r}^{s_{1}}\times\left(  B_{2,r}^{s_{2}%
}\right)  ^{n}$ with $\left(  \dot{B}_{2,1}^{\frac{n}{2}}\cap\dot{B}%
_{2,r}^{s_{1}}\right)  \times\left(  \dot{B}_{2,1}^{\frac{n}{2}}\cap\dot
{B}_{2,r}^{s_{2}}\right)  ^{n}$. For a definition and some basic properties
about homogeneous Besov spaces see \cite{Dan1}, Chapter $2$. To our knowledge,
$\left(  \dot{B}_{2,1}^{\frac{n}{2}}\cap\dot{B}_{2,r}^{s_{1}}\right)
\times\left(  \dot{B}_{2,1}^{\frac{n}{2}}\cap\dot{B}_{2,r}^{s_{2}}\right)
^{n}$ is the largest space in which one can prove long time existence for the
$abcd$ systems. Let us also point out that as opposed to classical Sobolev
spaces, working with Besov spaces enables us to attain the critical regularity
index $s=\frac{n}{2}+1$.

\subsection{The remaining cases\label{sub1}}

As we have seen earlier our method requires the restrictions \eqref{abcd}
imposed on the parameters $a,b,c,d$ in order to obtain the local existence
theory as well as estimate \eqref{I1}, see Remark \ref{obs}. However, the
general estimations \eqref{eq14} (see Section \ref{Section}) are valid for all
the parameters \eqref{1111} with $b+d>0$. Using \eqref{eq14} and supposing
that%
\begin{equation}
1+\varepsilon\eta_{0}\geq\alpha>0,\label{cond2}%
\end{equation}
we can obtain long time existence for any solutions of \eqref{eq2} in the
cases:
\begin{align*}
a &  =d=0,\text{ }b>0,\text{ }c<0\text{ and}\\
a &  =b=0,\text{ }d>0,\text{ }c<0.
\end{align*}
However the Friedrichs method is not well-suited in this case because when
establishing \eqref{eq14} we multiplied the frequency localized equation of
$V_{j}$ with $\eta V_{j}$. Thus we would run into trouble when establishing
the \eqref{eq14}-type estimate for the approximations $\left(  \eta^{m}%
,V^{m}\right)  $ introduced in Theorem \ref{Teorema1}, \eqref{En2}.
Nevertheless, one can imagine another strategy in order to bypass this
inconvenience namely establishing a convergence scheme for \eqref{eq2}. In
order to keep a certain homogeneity in this paper we will just give some brief
details of how this strategy would work. Let us take $W=0$ to fix the ideas.
Consider the linear systems:%
\begin{equation}
\left\{
\begin{array}
[c]{l}%
\left(  I-\varepsilon b\Delta\right)  \partial_{t}\eta^{m+1}%
+\operatorname{div}V^{m+1}+a\varepsilon\operatorname{div}\Delta V^{m+1}%
+\varepsilon\operatorname{div}\left(  \eta^{m}V^{m+1}\right)  =0,\\
\left(  I-\varepsilon d\Delta\right)  \partial_{t}V^{m+1}+\nabla\eta
^{m+1}+c\varepsilon\nabla\Delta\eta^{m+1}+\varepsilon\nabla V^{m+1}V^{m}=0,\\
\eta_{|t=0}^{m+1}=S_{m+1}\eta_{0}\text{ , }V_{|t=0}^{m+1}=S_{m+1}V_{0},
\end{array}
\right.  \tag{$\mathcal{L}_{m+1}$}%
\end{equation}
with $S_{m}=%
{\displaystyle\sum\limits_{j\leq m-1}}
\Delta_{j}$, and $\left(  \eta^{0},V^{0}\right)  =\left(  0,0\right)  $. We
claim that it is possible to establish the following \eqref{eq14}-type
estimate:
\begin{align*}
&  \frac{1}{2}\partial_{t}\left(  \int\left(  \eta^{m+1}\right)
^{2}+\varepsilon\left(  b-c\right)  \left\vert \nabla\eta^{m+1}\right\vert
^{2}+\varepsilon^{2}(-c)b\nabla^{2}\eta^{m+1}:\nabla^{2}\eta^{m+1}\right)  +\\
&  \frac{1}{2}\partial_{t}\left(  \int\left(  1+\varepsilon\eta^{m}\right)
\left\vert V^{m+1}\right\vert ^{2}+\varepsilon\left(  d-a+d\varepsilon\eta
^{m}\right)  \left\vert \nabla V^{m+1}\right\vert ^{2}+\varepsilon
^{2}(-a)d\nabla^{2}V^{m+1}:\nabla^{2}V^{m+1}\right)  \\
&  \leq\varepsilon P_{1}\left(  \left\Vert \left(  \eta^{m+1},V^{m+1}\right)
\right\Vert _{B_{2,r}^{s}},\sqrt{\varepsilon}\left\Vert \nabla\left(
\eta^{m+1},V^{m+1}\right)  \right\Vert _{B_{2,r}^{s}},\varepsilon\left\Vert
\nabla^{2}\left(  \eta^{m+1},V^{m+1}\right)  \right\Vert _{B_{2,r}^{s}%
}\right)  \times\\
&  P_{2}\left(  \left\Vert \left(  \eta^{m},V^{m}\right)  \right\Vert
_{B_{2,r}^{s}},\sqrt{\varepsilon}\left\Vert \nabla\left(  \eta^{m}%
,V^{m}\right)  \right\Vert _{B_{2,r}^{s}},\varepsilon\left\Vert \nabla
^{2}\left(  \eta^{m},V^{m}\right)  \right\Vert _{B_{2,r}^{s}}\right)  ,
\end{align*}
where $P_{1}$ and $P_{2}$ are two polynomials of degree $2$ with coefficients
depending only on the $a$, $b$, $c$, $d$ parameters, on the dimension $n$ and
on the regularity index $s$ but not on $\varepsilon$. Afterwards, Gronwall's
lemma combined with some stability estimates and a continuity argument would
allow us to obtain long time existence for \eqref{eq2} in the remaining two
cases. Of course, we would use in a decisive manner the non-cavitation
condition \eqref{cond2}.

\subsection{The $abcd$-system with smooth general topography\label{sub2}}

As we have mentioned in the introduction $abcd$-type models have been
established in \cite{Chaz} for a bottom given by the surface:%

\[
\{\left(  x,y,z\right)  :z=-h+\varepsilon S(x,y)\},
\]
where $S$ is smooth enough. In this case the $abcd$-system reads%
\begin{equation}
\left\{
\begin{array}
[c]{l}%
\left(  I-\varepsilon b\Delta\right)  \partial_{t}\eta+\operatorname{div}%
V+a\varepsilon\operatorname{div}\Delta V+\varepsilon\operatorname{div}\left(
(\eta-S)V\right)  =0,\\
\left(  I-\varepsilon d\Delta\right)  \partial_{t}V+\nabla\eta+c\varepsilon
\nabla\Delta\eta+\varepsilon\frac{1}{2}\nabla\left\vert V\right\vert ^{2}=0,\\
\eta_{|t=0}=\eta_{0},~V_{|t=0}=V_{0}.
\end{array}
\right.  \label{fin}%
\end{equation}
We claim that our method applies to this system with some minor modifications.
Let us give a few details. For the sake of simplicity we will suppose that
$\operatorname{curl}V_{0}=0$. By localizing the above equation in Fourier
space and proceeding in the same manner as in Section \ref{S2} we can see that
the terms preventing us from establishing long time existence are%
\[
\varepsilon\int\left(  \eta-S\right)  \eta_{j}\operatorname{div}V_{j}\text{
and }-c\varepsilon^{2}\int\left(  \eta-S\right)  \nabla\operatorname{div}%
V_{j}\nabla\eta_{j}.
\]
In order to repair this inconvenience we proceed as in Section \ref{refined},
the only difference being that we must multiply the localized equation of $V$
with $\left(  \eta-S\right)  V_{j}~$rather than simply $\eta V_{j}$. Thus, we
can obtain an estimate similar to \eqref{eq14} for the following quantity%
\begin{gather}
\frac{1}{2}\partial_{t}\left(  \int\eta_{j}^{2}+\varepsilon\left(  b-c\right)
\left\vert \nabla\eta_{j}\right\vert ^{2}+\varepsilon^{2}(-c)b\left(
\nabla^{2}\eta_{j}:\nabla^{2}\eta_{j}\right)  \right) \nonumber\\
+\frac{1}{2}\partial_{t}\left(  \int\left(  1+\varepsilon(\eta-S)\right)
V_{j}^{2}+\varepsilon\left(  d-a+d\varepsilon(\eta-S)\right)  \left(  \nabla
V_{j}:\nabla V_{j}\right)  +\varepsilon^{2}(-a)d\left(  \nabla^{2}V_{j}%
:\nabla^{2}V_{j}\right)  \right)  \leq\varepsilon G(U_{s})\nonumber
\end{gather}
where $G$ should be some polynomial function with its coefficients not
depending on $\varepsilon$. Afterwards proceeding like in the rest of the
paper, and taking $S\in B_{2,r}^{\tilde{s}}$ with some $\tilde{s}$ large
enough, one can obtain long time existence results for \eqref{fin} similar to
those obtained in Theorem \ref{Teorema2}.

\subsection*{Acknowledgement}

I would like to thank Prof. Rapha\"{e}l Danchin for the discussions and
suggestions that have improved the quality of this work. Also, I would like to
thank Prof. J.-C. Saut for having the kindness to present me the subjects of
interest along with his suggestions concerning the abcd-systems. \newpage

\section{Appendix: Littlewood-Paley theory}

We present here a few results of Fourier analysis used through the text. The
full proofs along with other complementary results can be found in
\cite{Dan1}. In the following if $\Omega\subset\mathbb{R}^{n}$ is a domain
then $\mathcal{D}\left(  \Omega\right)  $ will denote the set of smooth
functions on $\Omega$ with compact support and $\mathcal{S}$ will denote the
Schwartz class of functions defined on $\mathbb{R}^{n}$. Also, we consider
$\mathcal{S}^{\prime}$ is the set of tempered distributions on $\mathbb{R}%
^{n}$.

\begin{proposition}
\label{diadic}Let $\mathcal{C}$ be the annulus $\{\xi\in\mathbb{R}^{n}%
:3/4\leq\left\vert \xi\right\vert \leq8/3\}$. There exist two radial functions
$\chi\in\mathcal{D}(B(0,4/3))$ and $\varphi\in\mathcal{D(C)}$ valued in the
interval $\left[  0,1\right]  $ and such that:%
\begin{align}
\forall\xi &  \in\mathbb{R}^{n}\text{, \ }\chi(\xi)+\sum_{j\geq0}%
\varphi(2^{-j}\xi)=1\text{,}\label{25}\\
\forall\xi &  \in\mathbb{R}^{n}\backslash\{0\}\text{, \ }\sum_{j\in\mathbb{Z}%
}\varphi(2^{-j}\xi)=1\text{,}\label{26}\\
2  &  \leq\left\vert j-j^{\prime}\right\vert \Rightarrow\mathrm{Supp}%
(\varphi(2^{-j}\cdot))\cap\mathrm{Supp}(\varphi(2^{-j^{\prime}}\cdot
))=\emptyset\label{27}\\
j  &  \geq1\Rightarrow\mathrm{Supp}(\chi)\cap\mathrm{Supp}(\varphi(2^{-j}%
\cdot))=\emptyset\label{28}%
\end{align}
the set $\mathcal{\tilde{C}=B(}0,2/3\mathcal{)+C}$ is an annulus and we have%
\begin{equation}
\left\vert j-j^{\prime}\right\vert \geq5\Rightarrow2^{j}\mathcal{C}%
\cap2^{j^{\prime}}\mathcal{\tilde{C}}=\emptyset\text{.} \label{29}%
\end{equation}
Also the following inequalities hold true:%
\begin{align}
\forall\xi &  \in\mathbb{R}^{n}\text{, \ }\frac{1}{2}\leq\chi^{2}(\xi
)+\sum_{j\geq0}\varphi^{2}(2^{-j}\xi)\leq1\text{,}\label{210}\\
\forall\xi &  \in\mathbb{R}^{n}\backslash\{0\}\text{, \ }\frac{1}{2}\leq
\sum_{j\in\mathbb{Z}}\varphi^{2}(2^{-j}\xi)\leq1\text{.} \label{211}%
\end{align}

\end{proposition}

From now on we fix two functions $\chi$ and $\varphi$ satisfying the
assertions of the above proposition. The following two lemmas represent some
basic properties of the dyadic operators.

\begin{lemma}
For any $u\in\mathcal{S}^{\prime}$ we have:%
\[
u=\sum_{j\in\mathbb{Z}}\Delta_{j}u\text{ \ in \ }\mathcal{S}^{\prime}\left(
\mathbb{R}^{n}\right)  .
\]

\end{lemma}

\begin{lemma}
For any $u\in\mathcal{S}^{\prime}$ and $v\in\mathcal{S}^{\prime}$ we have
that:%
\[
\Delta_{j}\Delta_{\ell}u=0\text{ \ if \ }\left\vert j-\ell\right\vert \geq2.
\]

\end{lemma}

\subsection{Besov spaces}

The following propositions are a list of important basic properties of Besov
spaces that are used in the paper.

\begin{proposition}
\label{propA}A tempered distribution $u\in S^{\prime}$ belongs to $B_{2,r}%
^{s}$ if and only if there exists a sequence $\left(  c_{j}\right)  _{j} $
such that $\left(  2^{js}c_{j}\right)  _{j}\in\ell^{r}(\mathbb{Z)}$ with norm
$1$ and an universal constant $C>0$ such that for any $j\in\mathbb{Z}$ we have%
\[
\left\Vert \Delta_{j}u\right\Vert _{L^{2}}\leq Cc_{j}.
\]

\end{proposition}

\begin{proposition}
\label{PropBesov}Let $\ s,\tilde{s}\in\mathbb{R}$ and $r,\tilde{r}\in\left[
1,\infty\right]  $.

\begin{itemize}
\item $B_{2,r}^{s}$ is a Banach space which is continuously embedded in
$\mathcal{S}^{\prime}$.

\item The inclusion $B_{2,r}^{s}\subset B_{2,\tilde{r}}^{\tilde{s}}$ is
continuous whenever $\tilde{s}<s$ or $s=\tilde{s}$ and $\tilde{r}>r\,$.

\item We have the following continuous inclusion $B_{2,1}^{\frac{n}{2}}%
\subset\mathcal{C}_{0}$\footnote{$\mathcal{C}_{0}$ is the space of continuous
bounded functions which decay at infinity.}$(\subset L^{\infty})$.

\item (Fatou property) If $\left(  u_{n}\right)  _{n\in\mathbb{N}}$ is a
bounded sequence of $B_{2,r}^{s}$ which tends to $u$ in $\mathcal{S}^{\prime}$
then $u\in B_{2,r}^{s}$ and%
\[
\left\Vert u\right\Vert _{B_{2,r}^{s}}\leq\liminf_{n}\left\Vert u_{n}%
\right\Vert _{B_{2,r}^{s}}.
\]

\item If $r<\infty$ then
\[
\lim_{j\rightarrow\infty}\left\Vert u-S_{j}u\right\Vert _{B_{2,r}^{s}}=0.
\]

\end{itemize}
\end{proposition}

\begin{remark}
Taking advantage of the Fourier-Plancherel theorem and using \ref{211} one
sees that the classical Sobolev space $H^{s}$ coincides with $B_{2,2}^{s}$.
\end{remark}

\begin{proposition}
Let us consider $m\in\mathbb{R}$ and a smooth function $f:\mathbb{R}%
^{n}\rightarrow\mathbb{R}$ such that for all multi-index $\alpha$, there
exists a constant $C_{\alpha}$ such that:%
\[
\forall\xi\in\mathbb{R}^{n}\text{ \ }\left\vert \partial^{\alpha}f\left(
\xi\right)  \right\vert \leq C_{\alpha}\left(  1+\left\vert \xi\right\vert
\right)  ^{m-\left\vert \alpha\right\vert }.
\]
Then the operator $f\left(  D\right)  $ is continuous from $B_{2,r}^{s}$ to
$B_{2,r}^{s-m}$.
\end{proposition}

\begin{proposition}
\label{compact}Let $1\leq r\leq\infty$, $s\in\mathbb{R}$ and $\varepsilon>0$.
For all $\phi\in\mathcal{S}$, the map $u\rightarrow\phi u$ is compact from
$B_{2,r}^{s}$ to $B_{2,r}^{s-\varepsilon}$.
\end{proposition}

\begin{proposition}
\label{prodsmare} Let $s>0$ and $1\leq r\leq\infty$. Then $B_{2,r}^{s}\cap
L^{\infty}$ is an algebra. Moreover, we have:%
\[
\left\Vert uv\right\Vert _{B_{2,r}^{s}}\lesssim\left\Vert u\right\Vert
_{L^{\infty}}\left\Vert v\right\Vert _{B_{2,r}^{s}}+\left\Vert u\right\Vert
_{B_{2,r}^{s}}\left\Vert v\right\Vert _{L^{\infty}}.
\]
In particular, if $s>\frac{n}{2}$ or $s=\frac{n}{2}$ and $r=1$, $B_{2,r}^{s}$
is an algebra.
\end{proposition}

We end this section with the following result concerning a commutator-type
estimate. For a more general form of this result and its proof see \cite{Dan1}
page $116$, Lemma $2.100$.

\begin{proposition}
\label{comut} Let us consider $s\in\mathbb{R}$, $r\in\lbrack1,\infty]$ such
that $s>1+\frac{n}{2}$ or $s=1+\frac{n}{2}$ and $r=1$. Let $\left(
u,v\right)  \in B_{2,r}^{s}\times B_{2,r}^{s}$. We denote by%
\[
R_{j}=[\Delta_{j},u]\partial^{\alpha}v=\Delta_{j}\left(  u\partial^{\alpha
}v\right)  -u\Delta_{j}\partial^{\alpha}v,
\]
where $\alpha$ is any multi-index with $\left\vert \alpha\right\vert =1$.
Then, the following estimate holds true:%
\[
\left\Vert \left(  2^{js}\left\Vert R_{j}\right\Vert _{L^{2}}\right)
\right\Vert _{\ell^{r}(\mathbb{Z)}}\lesssim\left\Vert \nabla u\right\Vert
_{B_{2,r}^{s-1}}\left\Vert v\right\Vert _{B_{2,r}^{s}}.
\]

\end{proposition}


\newpage

\begin{thebibliography}{99}                                                                                               %


\bibitem {Amick}C.J. Amick, \emph{Regularity and uniqueness of solutions to
the Boussinesq system of equations}, J. Diff. Equations \textbf{54}, pp. 231-247.

\bibitem {Anh}C.T.Anh, \emph{On the Boussinesq/full dispersion and
Boussinesq/Boussinesq systems for internal waves}, Nonlinear Anal.,
\textbf{72} (2010), pp. 409-429.

\bibitem {Dan1}H. Bahouri, J.-Y. Chemin and R. Danchin, \emph{Fourier Analysis
and Nonlinear Partial Differential Equations}, Springer, Berlin 2010.

\bibitem {Bona1}J. L. Bona, M. Chen and J.-C. Saut, \emph{Boussinesq equations
and other systems for small-amplitude long waves in nonlinear dispersive
media. I: Derivation and linear theory}, J. Nonlinear Sci., \textbf{12
}(2002), 283-318.

\bibitem {Bona2}J. L. Bona, M. Chen and J.-C. Saut, \emph{Boussinesq equations
and other systems for small-amplitude long waves in nonlinear dispersive
media. II: The nonlinear theory}, Nonlinearity, \textbf{17 }(2004), 925-952.

\bibitem {Bona3}J. L. Bona, T. Colin, C. Guillop\'{e}, \emph{Propagation of
long-cested water waves}, Discrete and continuous dynamical systems,
\textbf{2} (2013), 599-628.

\bibitem {Bona4}J.L. Bona, T.Colin and D. Lannes, \emph{Long wave
approximation for water waves}, Archiv. Ration.Mech.Anal. \textbf{178} (2005), 373-410

\bibitem {Chen}M. Chen, \emph{Equations for bi-directional waves over an
uneven bottom}, Mathematical and Computers in Simulation \textbf{62} (2003),
pp. 3--9.

\bibitem {Chaz}F. Chazel, \emph{Influence of bottom topography on long water
waves}, M2AN Math. Model. Numer. Anal. \textbf{41} (4) (2007), 771--799.

\bibitem {Saut2}F. Linares, D. Pilod and J.-C. Saut, \emph{Well-Posedness of
strongly dispersive two-dimensional surface wave Boussinesq systems}, SIAM J.
Math. Anal., \textbf{44}, No. 6(2012) pp. 4195-4221

\bibitem {Lannes}D. Lannes, \emph{The Water Waves Problem: Mathematical
Analysis and Asymptotics}, vol. 188 of Mathematical Surveys and Monographs,
AMS, 2013.

\bibitem {Saut1}J.-C. Saut, Li Xu, \emph{The Cauchy problem on large time for
surface waves Boussinesq systems, }J. Math. Pures Appl. 97 (2012) 635--662.

\bibitem {Schon}M.E. Schonbek, \emph{Existence of solutions for the Boussinesq
system of equations}, J. Diff. Equations, \textbf{42}, pp.325-352.
\end{thebibliography}
\end{document}